\documentclass[11pt, twoside]{amsart}

\usepackage[utf8]{inputenc}
\usepackage[T1]{fontenc}
\usepackage{amssymb,amsmath,amstext}
\usepackage{hyperref, graphicx}

\newtheorem{theor}{Theorem}[section]
\newtheorem{lem}[theor]{Lemma}
\newtheorem{defin}[theor]{Definition}

\newtheorem{prop}[theor]{Proposition} 

\newtheorem{notation}[theor]{Notation}
\newtheorem{exam}[theor]{Example}

\newtheorem{cor}[theor]{Corollary}
\newtheorem{rem}[theor]{Remark}

\newtheorem{fact}[theor]{Fact}
\newtheorem{claim}[theor]{Claim}

\newtheorem{assump}[theor]{Assumption}

\numberwithin{equation}{section}

\newcommand{\acl}{\mathrm{acl}}

\newcommand{\tp}{\mathrm{tp}}

\newcommand{\cl}{\mathrm{cl}}

\newcommand{\es}{\emptyset}

\newcommand{\rk}{\mathrm{rk}}

\newcommand{\meq}{^{\mathrm{eq}}}

\newcommand{\mcB}{\mathcal{B}}

\newcommand{\mcM}{\mathcal{M}}

\newcommand{\mbbN}{\mathbb{N}}

\newcommand{\ind}{\raisebox{-2pt}[5pt][0pt]{$\smile$} \hspace*{-6.8pt}\raisebox{3pt}[5pt][0pt]{$|$} \; \: }
\newcommand{\nind}{\raisebox{-2pt}[5pt][0pt]{$\smile$} 
\hspace*{-6.8pt}\raisebox{3pt}[5pt][0pt]{$|$}\hspace*{-6.8pt}
\raisebox{3pt}[5pt][0pt]{$\diagup$} }

\newcommand{\indn}{\raisebox{-2pt}[5pt][0pt]{$\smile$} \hspace*{-6.8pt}\raisebox{3pt}[5pt][0pt]{$|^n$} }

\newcommand{\nindn}{\raisebox{-2pt}[5pt][0pt]{$\smile$} 
\hspace*{-6.8pt}\raisebox{3pt}[5pt][0pt]{$|^n$}\hspace*{-12pt}
\raisebox{3pt}[5pt][0pt]{$\diagup$} }

\newcommand{\indo}{\raisebox{-2pt}[5pt][0pt]{$\smile$} \hspace*{-6.8pt}\raisebox{3pt}[5pt][0pt]{$|^0$} }

\newcommand{\nindo}{\raisebox{-2pt}[5pt][0pt]{$\smile$} 
\hspace*{-6.8pt}\raisebox{3pt}[5pt][0pt]{$|^0$}\hspace*{-12pt}
\raisebox{3pt}[5pt][0pt]{$\diagup$} }

\newcommand{\nindone}{\raisebox{-2pt}[5pt][0pt]{$\smile$} 
\hspace*{-6.8pt}\raisebox{3pt}[5pt][0pt]{$|^1$}\hspace*{-12pt}
\raisebox{3pt}[5pt][0pt]{$\diagup$} }

\newcommand{\thind}{\raisebox{-2pt}[5pt][0pt]{$\smile$} 
\hspace*{-6.8pt}\raisebox{3pt}[5pt][0pt]{$|^{\text{\rm \th}}$}  }

\newcommand{\thnind}{\raisebox{-2pt}[5pt][0pt]{$\smile$} 
\hspace*{-6.8pt}\raisebox{3pt}[5pt][0pt]{$|^{\text{\rm \th}}$}\hspace*{-12pt}
\raisebox{3pt}[5pt][0pt]{$\diagup$} }

\newcommand{\rng}{\mathrm{rng}}

\title[Notions of rank and independence]
{Notions of rank and independence in countably categorical theories}

\author{Vera Koponen}

\address{Vera Koponen, Department of Mathematics, Uppsala University, Sweden.}

\email{vera.koponen@math.uu.se}

\date{25 May, 2026}

\begin{document}

\maketitle

\begin{abstract}
For an $\omega$-categorical theory $T$ and model $\mathcal{M}$ of $T$
we define a hierarchy of ranks, the $n$-ranks for $n < \omega$
which only care about imaginary elements ``up to level $n$'',
where level $n$ contains every element of $M$ and every imaginary element
that is an equivalence class of an $\es$-definable equivalence relation on $n$-tuples of elements from $M$.
Using the $n$-rank we define the notion of $n$-independence.
For all $n < \omega$, the $n$-independence relation restricted to $M_n$ 
has all properties of an independence relation according
to Kim and Pillay \cite{KP} with the {\em possible exception} of the symmetry property.
We prove that, given any $n < \omega$,
if $\mathcal{M} \models T$ and the algebraic closure in $\mathcal{M}^{\text{\rm eq}}$
restricted to imaginary elements ``up to level $n$'' 
which have $n$-rank 1 (over some set of parameters) satisfies the exchange property,
then $n$-independence is symmetric and hence an independence relation when restricted to $M_n$.
Then we show that if $n$-independence is symmetric for all $n < \omega$, then $T$ is rosy.
An application of this is that if $T$ has geometric elimination of imaginaries and the algebraic closure in $\mathcal{M}$
restricted to elements of $M$ of 0-rank 1 (over some set of parameters from $M^{\text{\rm eq}}$) satisfies the 
exchange property, then $T$ is superrosy with finite $U^{\text{\rm \th}}$-rank.
\end{abstract}

\section{Introduction}

\noindent
A variety of notions of rank and independence have played
an important role in model theory at least since Morley's influential work on uncountably categorical theories
\cite{Mor} in the 1960ies.
Such notions have been central for developing more or less general (meta) theories which divide complete
first-order theories into various classes. Shelah \cite{She} used them in his development of stability theory 
by which theories can be classified into $\omega$-stable, superstable, stable, or unstable.
Later Shelah's classification theory, and its notion of independence, was shown,
by Kim and Pillay \cite{KP, Kim_book}, to make sense in a wider context
and the classes of simple and supersimple first-order theories were introduced.
Yet later, more general notions of independence, including thorn-independence, have been studied
by Onshuus, Ealy and Adler
and the classes of rosy and superrosy theories have been introduced \cite{Adl, EO, Ons}.
The class of rosy theories is quite diverse and includes, for example, 
all simple theories (which exclude a linearly ordered universe) and all o-minimal theories 
(which assume a linearly ordered universe) \cite{Van}.
It is also the largest class of theories for which there is an independence relation which satisfies certain
basic and natural properties \cite{EO}.
The work on NIP theories \cite{Sim_book}
also uses notions of rank and independence (and NIP theories and rosy theories partially overlap).

Another direction of model theoretic research, which has used notions of rank and
independence as a crucial tool, 
has focused on understanding the {\em fine structure} of the models of more specific theories,
and on finding {\em ``nice'' axiomatizations} of such theories.
This line of research includes work on
totally categorical theories, uncountably categorical theories \cite{Zil, AZ, Hru1, Hru2},
theories of stable finitely homogeneous structures \cite{Lach97},
$\omega$-categorical $\omega$-stable theories \cite{CHL},
theories of smoothly approximable structures \cite{CH, KLM},
theories of simple finitely homogeneous structures \cite{BFM, Kop},
theories of Fra\"{i}ss\'{e} limits of classes of finite structures with the free amalgamation property
\cite{Con},
$\omega$-categorical NIP theories \cite{Sim22}, 
and NIP finitely homogeneous rosy theories \cite{OS}.
In the present context it is relevant that all
finitely homogeneous structures (which are also called ultrahomogeneous, or simply homogeneous,
and which include all Fra\"{i}ss\'{e} limits with a finite relational vocabulary) 
and all smoothly approximable structures are $\omega$-categorical.
 
In this study we consider (only) $\omega$-categorical theories $T$ and investigate a hierarchy of ranks,
the {\em $n$-ranks} for $n < \omega$, where the $n$-rank is defined entirely in terms of the algebraic closure
operator on $\mcM\meq$ (where $\mcM$ is a model of the theory) 
{\em restricted} to the set $M_n$ containing every real element and every imaginary element
that corresponds to an equivalence class of an $\es$-definable equivalence relation on 
$M^k = \underbrace{M \times \ldots \times M}_{\text{$k$ times}}$ where $k \leq n$.
With the $n$-rank we define a notion of {\em $n$-independence}, denoted $\indn$.
Without further assumptions than $\omega$-categoricity, $n$-independence restricted to $M_n$
has all the properties of an independence relation according to \cite{KP}
(possibly)
{\em except} for the symmetry property.
The main technical contribution of this study is to isolate a property, parametrized by $n$, which implies
(in fact, is equivalent to) 
that $n$-independence is symmetric.
The property in question 
(Assumption~\ref{assumption about exchange property})
is the {\em exchange property} of algebraic closure restricted to
elements of $M_n$ with $n$-rank 1 (over some $C \subseteq M\meq$).
More precisely, if, for any $C \subseteq M\meq$,
the algebraic closure using parameters from $C$ restricted to elements of $M_n$ with $n$-rank 1 over $C$
has the exchange property, then $n$-independence is symmetric, and hence it is an independence relation
in the sense of \cite{KP} when restricted to $M_n$, which is stated by
Theorem~\ref{n-independence is an independence relation}.
Using this we can relatively easily show the statement
of Theorem~\ref{the assumption implies that T is rosy}
that if, for {\em all} $n < \omega$, the algebraic closure restricted to elements of $M_n$ with $n$-rank 1
has the exchange property,
then the theory is rosy, which is proved by showing that thorn-independence has local character
(see \cite[Theorem~3.7]{EO}).
Thus, an $\omega$-categorical theory which is {\em not} rosy must, for some $n < \omega$,
have a model $\mcM$ and $C \subseteq M\meq$ such that 
algebraic closure using parameters from $C$ restricted to elements of $M_n$ with $n$-rank 1 over $C$
does {\em not} have the exchange property.
This may be useful for finding a dividing line between $\omega$-categorical rosy theories 
(e.g. dense linear order and the Fr\"{i}ss\'{e} limit of a class of finite relational structures 
with free amalgamation \cite{Con}) and 
$\omega$-categorical non-rosy theories (e.g. the Fra\"{i}ss\'{e} limit of finite boolean algebras and
$T^*_{feq}$ in \cite[Example 2.6]{Adl}).
A consequence of 
Theorem~\ref{the assumption implies that T is rosy}
is that if algebraic closure on $M\meq$ is trivial in the sense of 
Definition~\ref{definition of trivial acl up to n},
then the theory is rosy, as stated by
Theorem~\ref{trivial acl implies that T is rosy}.
 
We then consider the effects of assuming geometric elimination of imaginaries
(Definition~\ref{definition of geometric elimination of hyperimaginaries})
in the context studied here,
where geometric elimination of imaginaries is a consequence of the perhaps more familiar concept of
weak elimination of imaginaries \cite{HHM, Hod, NW, Con}.
It will be shown that if the $\omega$-categorical theory has geometric elimination of imaginaries, 
then the hierarchy of $n$-independence notions, for $n < \omega$, collapses to the bottom level,
that is for all $n < \omega$, $A, B, C \subseteq M\meq$, $A \underset{C}{\indn} B$
if and only if $A \underset{C}{\indo} B$.

We use this to show that if the $\omega$-categorical theory $T$
has geometric elimination of imaginaries
and, for any $\mcM \models T$ and $C \subseteq M\meq$,
algebraic closure using parameters from $C$ restricted to elements of $M$ (i.e ``real'' elements)
that have 0-rank 1 over $C$
has the exchange property,
then $T$ is superrosy with finite $U^{\text{\rm \th}}$-rank;
this is Theorem~\ref{concluding that T is superrosy} below.
Hence Theorems~\ref{the assumption implies that T is rosy} and~\ref{concluding that T is superrosy}
are generalizations, in the context of $\omega$-categorical theories, 
of the result that if a theory has geometric elimination of imaginaries and algebraic closure on 
all real elements (not just those of $0$-rank 1 over some set) has the exchange property,
then the theory is superrosy with $U^{\text{\rm \th}}$-rank 1
(this is \cite[Theorem~4.12]{EO} where the authors also contribute the result to
Gagelman \cite{Gag} and Adler).

Conant \cite{Con} has proved that the theory of the Fra\"{i}ss\'{e} limit
of any class of finite (relational) structures with {\em free} amalgamation has weak elimination of imaginaries and 
is superrosy with  $U^{\text{\rm \th}}$-rank 1. 
Such a theory has trivial algebraic closure on real elements and therefore 
the algebraic closure on real elements satisfies the exchange property.
So unless every theory to which 
Theorem~\ref{concluding that T is superrosy}
applies can be constructed as the Fra\"{i}ss\'{e} limit
of  a class of relational structures with {\em free} amalgamation, 
then Theorem~\ref{concluding that T is superrosy} applies to a larger class of theories than \cite{Con}.
It may also be the case that Theorem~\ref{the assumption implies that T is rosy}
or Theorem~\ref{concluding that T is superrosy} can be applied
to $\omega$-categorical theories obtained by Hrushovski's method of construction of a ``generic structure''
for a class of finite structures
(see e.g. \cite{Eva}),
but I have not investigated this.

\section{Preliminaries}

\noindent
We assume familiarity with basic model theory as can be found in for example \cite{Hod, She, TZ}.
A first-order structure will be denoted by $\mcM$ and its universe (also called domain) by $M$.
Subsets of the universe (of a first-order structure) 
will be denoted by $A, B, C, D$ and {\em finite} sequences (also called tuples)
of elements of the universe will be denoted by $\bar{a}, \bar{b}, \bar{c}$, etcetera.
If $S$ is a set then $|S|$ denotes its cardinality, and if $\bar{s}$ is a finite sequence then $|\bar{s}|$ denotes its length
and $\rng(\bar{s})$ denotes the set of elements occuring in $\bar{s}$.
By $Th(\mcM)$ we denote the complete first-order theory of the structure $\mcM$.
We sometimes write $AB$ to denote the union $A \cup B$ of the sets $A$ and $B$, 
and if $\bar{a} = (a_1, \ldots, a_k)$ we sometimes write $\bar{a}B$ to denote the
set
$\{a_1, \ldots, a_k\} \cup B$.

We assume familiarity with the structure $\mcM\meq$ ``with imaginaries'' which is constructed from
any structure $\mcM$ as explained in (for example) \cite{She} and \cite{Hod}.
In this study we will be concerned with ``$n$-level approximations'' of $M\meq$ (the universe of $\mcM\meq$)
as follows:

\begin{defin}\label{definition of M-n}{\rm
Let $V$ be a vocabulary and let $\mcM$ be a $V$-structure with universe $M$.
Let $V \cup W$ be the vocabulary of $\mcM\meq$.
Define $M_0 = M$ and, for all $n < \omega$,
$M_{n+1} = M_n \cup X_{n+1}$ where $X_{n+1}$ is the set of all imaginary elements $a \in M\meq$
such that $a$ is an equivalence class of a ($\es$-definable) equivalence relation on 
$M^{n+1} = \underbrace{M \times \ldots \times M}_{\text{$n+1$ times}}$.
For each $n < \omega$, let $\mcM'_n$ be the reduct of $\mcM\meq$ to the subvocabulary of $V \cup W$
that contains $V$ and every symbol from $W$ that is associated to a $\es$-definable equivalence relation on $M^k$
for some $1 \leq k \leq n$. Then let $\mcM_n$ be the substructure of $\mcM'_n$ with universe $M_n$.
}\end{defin}

\begin{defin}\label{definition of acl-n}{\rm
Let $\mcM$ be a structure.
For every $A \subseteq M\meq$, $\acl\meq(A)$ denotes the algebraic closure of $A$ computed in $\mcM\meq$.
For all $n < \omega$ and $A \subseteq M\meq$, we define
$\acl^n(A) = \acl\meq(A) \cap M_n$.
}\end{defin}

\noindent
The following follows directly from the definition above:

\begin{lem}\label{acl n-0 and acl n}
Suppose that $m < n < \omega$ and $C \subseteq M\meq$.\\
(i)  If $a \in \acl^m(C)$ then $a \in \acl^n(C)$.\\
(ii) If $a \in M_m$ and $a \notin \acl^m(C)$ then $a \notin \acl^n(C)$.
\end{lem}

\noindent
Since each $\mcM_n$ is interpretable in $\mcM$ the following result follows from \cite[Theorem~7.3.8]{Hod}:

\begin{fact}\label{M-n is omega-categorical}
If $T$ is $\omega$-categorical and $\mcM \models T$ then, for every $n < \omega$, 
$Th(\mcM_n)$ is $\omega$-categorical.
\end{fact}

\noindent 
We will use the following facts about $\omega$-categorical theories (briefly explained below):

\begin{fact}\label{consequences of omega-categoricity}
(i) If $T$ is $\omega$-categorical then, for all $0 < n < \omega$, 
there are only finitely many complete types over $\es$ in the free variables 
$x_1, \ldots, x_n$ and each one of them is implied (modulo $T$) by a single formula (the formula that isolates the type).
Hence, there are only finitely many nonequivalent (modulo $T$) formulas in the same free variables.\\
(ii) If $T$ is $\omega$-categorical, $\mcM \models T$, $n < \omega$, and $B \subseteq M\meq$ is finite,
then $\acl^n(B)$ is finite.
\end{fact}

\noindent
The first part above follows from the well-known theorem of
Engeler, Ryll-Nardzewski, and Svenonius \cite[Theorem~7.3.1]{Hod}.
The second part is a consequence of the fact that, since $B$ is finite, there is $n \leq m < \omega$ such that
$B \subseteq M_m$ and $Th(\mcM_m)$ is $\omega$-categorical.

\subsection*{In the rest of the article we make the following assumptions}

\begin{enumerate}
\item $T$ is a complete $\omega$-categorical theory in a countable language.
\item $\mcM \models T$ is $\kappa$-saturated where $\kappa$ is an infinite cardinal (which can be chosen as large as
we like).
It follows that $\mcM\meq$ is $\kappa$-saturated and we let $T\meq$ be the complete theory of $M\meq$.
\item All subsets of $M\meq$ that are mentioned have cardinality less than $\kappa$.
\end{enumerate}

\begin{notation}{\rm
If $a_1, \ldots, a_n \in M\meq$ and $B  \subseteq M\meq$, then $\tp(a_1, \ldots, a_n / B)$
denotes the complete type of $a_1, \ldots, a_n$ over $B$ computed in $\mcM\meq$, 
and $\tp(a_1, \ldots, a_n)$ is an abbreviation of $\tp(a_1, \ldots, a_n / \es)$.
}\end{notation}

\begin{defin}\label{definition of independence relation}{\rm
(i) Let us call a {\em bijective} function $\sigma$ from a subset $A \subseteq M\meq$ to a subset of $M\meq$
{\em elementary} if for all $n < \omega$ and all $a_1, \ldots, a_n \in A$,
$\tp(a_1, \ldots, a_n) = \tp(\sigma(a_1), \ldots, \sigma(a_n))$.\\
(ii) Following \cite{KP} we say that a collection $\Gamma$ of triples $(A, B, C)$, 
where 
$A, B, C \subseteq M\meq$ and {\em $A$ is finite}, 
is an {\em independence relation} if the following hold, where 
$A \underset{C}{\ind} B$ means that $(A, B, C) \in \Gamma$:
\begin{enumerate}
\item (invariance) If $\sigma$ is an elementary function from some subset of $M\meq$
that includes $ABC$ to some subset of $M\meq$,
then $A \underset{C}{\ind} B$ if and only if $\sigma(A) \underset{\sigma(C)}{\ind} \sigma(B)$.

\item (local character) For all $A$ and $B$ (where $A$ is finite) there is a countable $C \subseteq B$ such that
$A \underset{C}{\ind} B$.

\item (finite character) $A \underset{C}{\ind} B$ if and only if for all finite $B' \subseteq B$,
$A \underset{C}{\ind} B'$.

\item (extension) For all $n < \omega$, 
$A = \{a_1, \ldots, a_n\}$, $B$, and $C$, such that $C \subseteq B$, there is $A' = \{a'_1, \ldots, a'_n\}$ such that
$\tp(a'_1, \ldots, a'_n / C) = \tp(a_1, \ldots, a_n / C)$ and $A' \underset{C}{\ind} B$.

\item (monotonicity)\footnote{
In the formulation of \cite{KP} this condition is a part of the transitivity property.}
 If $B \subseteq C \subseteq D$ and
$A \underset{B}{\ind} D$ then $A \underset{B}{\ind} C$ and $A \underset{C}{\ind} D$.

\item (transitivity) If $B \subseteq C \subseteq D$,
 $A \underset{B}{\ind} C$ and $A \underset{C}{\ind} D$, then $A \underset{B}{\ind} D$.

\item (symmetry) For all finite $A$ and $B$ and any $C$, $A \underset{C}{\ind} B$ if and only if
$B \underset{C}{\ind} A$.
\end{enumerate}
(iii) We call $\Gamma$ an {\em independence relation restricted to $M_n$} 
if (1) -- (7) hold whenever $A, B, C, D \subseteq M_n$.
}\end{defin}

\begin{rem}\label{remark on locality}{\rm
The local character of $\ind$ (property~(2) above) is (assuming that $T$ is countable) equivalent to the following:
There is do not exist finite $A \subseteq M\meq$ and finite $B_\alpha$,
for all $\alpha < \aleph_1$, such that $A \underset{\bigcup_{i < \alpha}B_i}{\nind} B_\alpha$ for all $\alpha < \aleph_1$.
This formulation (roughly) is used in e.g. \cite{EO}, but I have not found a clear statement in the literature of the
equivalence of the two versions of local character. 
However, the argument in the proof of Proposition~2.3.7 in \cite{Kim_book} (the part showing that
conditions~(1) and~(2) of that proposition are equivalent) can easily be adapted to prove the equivalence
of the two versions of local character.
}\end{rem}

\begin{defin}\label{definition of geometric elimination of hyperimaginaries}{\rm
We say that $T$ has {\em geometric elimination of imaginaries} if for all $a \in M\meq$ there is
a finite sequence $\bar{b}$ of elements from $M$ such that 
$a \in \acl\meq(\bar{b})$ and $\rng(\bar{b}) \subseteq \acl\meq(a)$, or equivalently,
$\acl\meq(a) = \acl\meq(\bar{b})$.
}\end{defin}

\noindent
The more commonly used notion of {\em weak elimination of imaginaries} 
(where, with the notation of the definition above, 
it is required that $a$ belongs to the definable closure of $\bar{b}$) implies geometric elimination of imaginaries.
Weak elimination of imaginaries has been studied, for example, 
by Hodges, Hodkinson, and Macpherson, in \cite{HHM} where they demonstrated,
among other things, that
the complete theories of dense linear order, the random (or Rado) graph, and the Fra\"{i}ss\'{e} limit of the set
of finite $K_n$-free graphs have weak elimination of imaginaries.
Conditions that apply to $\omega$-categorical theories and 
imply weak elimination of imaginaries are given by \cite[Lemma~6.5 and Proposition~8.2]{HHM}.
Newelski and Wencel \cite{NW} have proved that the complete theory of an infinite boolean algebra 
with only finitely many atoms has weak elimination of imaginaries.
The above mentioned result about $K_n$-free graphs has been generalized by
Conant \cite{Con} to a result saying that the theory of every Fra\"{i}ss\'{e} limit  of a class of
finite structures with free amalgamation has weak elimination of imaginaries.

The backbone of the theory that will be developed does not use the assumption that $T$ has geometric elimination of imaginaries,
but, as we will see, the ``hierarchy'' of ranks and independence relations that we will study 
collapses to the bottom level under
the assumption of geometric elimination of imaginaries.

In Section~\ref{The exchange property and symmetry}
we will need the following concept \cite{Hod, Ox, TZ}:

\begin{defin}\label{definition of pregeometry}{\rm 
A {\em pregeometry}, also called {\em matroid}, consists of a set $X$ and a function
$\cl$ from the powerset of $X$ to the powerset of $X$ which has the following properties:
\begin{enumerate}
\item For all $A \subseteq X$, $A \subseteq \cl(A)$.
\item For all $A \subseteq X$, $\cl(\cl(A)) = \cl(A)$.
\item For all $A \subseteq X$, if $a \in \cl(A)$ then $a \in \cl(A')$ for some {\em finite} $A' \subseteq A$.
\item (Exchange property) For all $A \subseteq X$ and all $b, c \in X$, if $b \in \cl(A \cup \{c\}) \setminus \cl(A)$,
then $c \in \cl(A \cup \{b\})$.
\end{enumerate}
}\end{defin}

\begin{fact}\label{basis of a pregeometry}
Let $(X, \cl)$ be a pregeometry and let $A \subseteq X$.\\
(i) Then there is $B \subseteq A$ such that $A \subseteq \cl(B)$ 
(``$B$ spans $A$'') and for every
$b \in B$, $b \notin \cl(B \setminus \{b\})$ (``$B$ is independent'').
We call such $B$ a {\em basis} of $A$.\\
(ii) All bases of $A$ have the same cardinality which we call the {\em dimension} of $A$.
\end{fact}

\section{$n$-rank}\label{properties of n-rank}

\noindent
In this section we develop, for an arbitrary $n < \omega$, a theory of a notion of ``$n$-rank''
which will be used, in the next section, to define a notion of ``$n$-independence''.

\begin{lem}\label{bound on length of special sequences}
Let $n < \omega$, let $A  \subseteq M_n$ be finite and let $B \subseteq M\meq$.
Then there is $r < \omega$, depending only on $n$ and $|A|$, such that
if $a_1, \ldots, a_m \in \acl^n(A)$ and, for all $k = 1, \ldots, m$,
$a_k \notin \acl\meq(a_1, \ldots, a_{k-1}, B)$, 
then $m \leq r$.
\end{lem}

\noindent
{\bf Proof.}
Let $n < \omega$, let $A  \subseteq M_n$ be finite, with cardinality $s < \omega$, say, and let $B \subseteq M\meq$.
Suppose that $a_1, \ldots, a_m \in \acl^n(A)$ and, for all $k = 1, \ldots, m$,
$a_k \notin \acl\meq(a_1, \ldots, a_{k-1},  B)$.
Then $a_i \neq a_j$ if $i \neq j$ so $|\acl^n(A)| \geq m$.
Since $Th(\mcM_n)$ is $\omega$-categorical 
(by Fact~\ref{M-n is omega-categorical}) 
it follows there is $r < \omega$ such that whenever $A' \subseteq M^n$ and
$|A'| \leq s$, then $|\acl^n(A')| \leq r$. 
Hence $m \leq |\acl^n(A)| \leq r$.
\hfill $\square$

\begin{defin}\label{definition of n-rank}{\rm
Let $n < \omega$.
For all $A, B \subseteq M\meq$ we define the {\em $n$-rank of $A$ over $B$}, denoted $\rk^n(A / B)$, as follows:
\begin{enumerate}
\item $\rk^n(A / B) \geq 0$.
\item For any ordinal $\alpha$, $\rk^n(A / B) \geq \alpha + 1$ if there is $a \in \acl^n(A) \setminus \acl^n(B)$ such that
$\rk^n(A / \{a\} \cup B) \geq \alpha$.
\item For a limit ordinal $\alpha$, $\rk^n(A / B) \geq \alpha$ if $\rk^n(A / B) \geq \beta$ for all $\beta < \alpha$.
\end{enumerate}
Finally, $\rk^n(A / B) = \alpha$ if $\rk^n(A / B) \geq \alpha$ and $\rk^n(A / B) \not\geq \alpha + 1$.
If $\rk^n(A / B) \geq \alpha$ for all ordinals $\alpha$ we say that the $n$-rank of $A$ over $B$ is {\em undefined}.
We define $\rk^n(A) = \rk^n(A / \es)$.
If $\bar{a}$ is a sequence of elements (from $M\meq$) then $\rk^n(\bar{a} / B)$ means the same as $\rk^n(\rng(\bar{a}) / B)$,
and if $\bar{b}$ is a sequence of elements then $\rk^n(\bar{a} / \bar{b})$ means the same as $\rk^n(\rng(\bar{a}) / \rng(\bar{b}))$.
}\end{defin}

\noindent
Note that it follows that $\rk^n(A / B) \geq 1$ if and only if $\acl^n(A) \not\subseteq \acl^n(B)$.

\medskip
\noindent
{\em For the rest of this section we fix an arbitrary $n < \omega$.}

\begin{exam}\label{example for n-rank}{\rm
Let $T$ be the theory which expresses that $E$ is an equivalence relation with infinitely many equivalence
classes all of which are infinite.
Let $\mcM \models T$ and let $a \in M$.
It is well-known (and easy to show) that $T$ is $\omega$-categorical with elimination of quantifiers.
From this it easily follows that, for all $A \subseteq M$, $\acl^0(A) = A$.
Therefore $\rk^0(a) = 1$.
Let $[a]_E$ denote the equivalence class of $a$ with respect to $E$ as a element of $M\meq$, so $[a]_E \in M_1$.
Then $\rk^1(a / \{a, [a]_E\}) \geq 0$, and as $a \in \acl^1(a) \setminus \acl^1(\es)$ we 
get $\rk^1(a / \{[a]_E\}) \geq 1$.
Since $[a]_E \in \acl^1(a) \setminus \acl^1(\es)$ we get $\rk^1(a) \geq 2$.
As $T$ has elimination of quantifiers it follows that $E$ is the only $\es$-definable equivalence relation on $M$.
Hence $\acl^1(a) = \{a, [a]_E\}$ and consequently $\rk^1(a) = 2$.
Let $b \in M$ be such that $b \neq a$ and $[b]_E = [a]_E$.
By arguing similarly as above it follows that $\rk^0(a / b) = \rk^0(a) = 1$
and $\rk^1(a / b) = 1 < 2 = \rk^1(a)$.
}\end{exam}

\begin{lem}\label{n+1-rank is at least as big as n-rank}
Let $A, B \subseteq M\meq$ and suppose that $\rk^n(A / B)$ and $\rk^{n+1}(A / B)$ are defined.
Then  $\rk^n(A / B) \leq \rk^{n+1}(A / B)$.
\end{lem}

\noindent
{\bf Proof.}
We prove by induction that if $\rk^n(A / B) \geq \alpha$ then $\rk^{n+1}(A / B) \geq \alpha$.
This is clear for $\alpha = 0$.
So suppose that $\rk^n(A / B) \geq \alpha + 1$. 
Then there is $a \in \acl^n(A) \setminus \acl^n(B)$ such that $\rk^n(A / aB) \geq \alpha$.
By the induction hypothesis we get $\rk^{n+1}(A / aB) \geq \alpha$.
We have $a \in \acl^n(A) \subseteq \acl^{n+1}(A)$
and $a \notin \acl^n(B)$ so (by Lemma~\ref{acl n-0 and acl n}) $a \notin \acl^{n+1}(B)$.
Hence $\rk^{n+1}(A / B) \geq \alpha + 1$.
If $\alpha$ is a limit ordinal and $\rk^n(A / B) \geq \alpha$ then $\rk^n(A / B) \geq \beta$ for all $\beta < \alpha$,
so by the induction hypothesis $\rk^{n+1}(A / B) \geq \beta$ for all $\beta < \alpha$, hence
$\rk^{n+1}(A / B) \geq \alpha$.
\hfill $\square$

\begin{lem}\label{characterization of rank}
Let $\alpha < \omega$ and $A, B \subseteq M\meq$. \\
(i) $\rk^n(A / B) \geq \alpha$ if and only if there are
$a_1, \ldots, a_\alpha \in \acl^n(A)$ such that, for all $k = 1, \ldots, \alpha$,
$a_k \notin \acl^n(\{a_1, \ldots, a_{k-1}\} \cup B)$. \\
(ii) If $\rk^n(A / B) = \alpha$ then $\alpha$ is maximal such that there are
$a_1, \ldots, a_\alpha \in \acl^n(A)$ such that, for all $k = 1, \ldots, \alpha$,
$a_k \notin \acl^n(\{a_1, \ldots, a_{k-1}\} \cup B)$.\\
(iii) Suppose that $\rk^n(A / B) = \alpha$, $a_1, \ldots, a_\alpha \in \acl^n(A)$ and, for all $k = 1, \ldots, \alpha$,
$a_k \notin \acl^n(\{a_1, \ldots, a_{k-1}\} \cup B)$.
Then, for all $k = 1, \ldots, \alpha$,
\begin{enumerate}
\item $\rk^n(a_k / \{a_1, \ldots, a_{k-1}\} \cup B) = 1$,
\item $\rk^n(A / \{a_1, \ldots, a_k\} \cup B) = \alpha - k$, and
\item $\rk^n(a_1, \ldots, a_k / B) = k$.
\end{enumerate}
In particular we have $\rk^n(A / \{a_1, \ldots, a_\alpha\} \cup B) = 0$ so
$\acl^n(A) \subseteq \acl^n(\{a_1, \ldots, a_\alpha\} \cup B)$, and
$\rk^n(a_1, \ldots, a_\alpha / B) = \alpha = \rk^n(A / B)$.
\end{lem}

\noindent
{\bf Proof.}
Let $\alpha < \omega$ and $A, B \subseteq M\meq$. 
We prove (i) by induction on $\alpha$.
For $\alpha = 0$ the statement is vacuous.
Suppose that $\rk^n(A / B) \geq \alpha + 1$.
Then there is $a \in \acl^n(A) \setminus \acl^n(B)$ such that $\rk^n(A / aB) \geq \alpha$.
By the induction hypothesis there are $a_1, \ldots, a_\alpha \in \acl^n(A)$ such that, for all  $k = 1, \ldots, \alpha$,
$a_k \notin \acl^n(\{a, a_1, \ldots, a_{k-1}\} \cup B)$.
If we rename $a_i$ by $a_{i+1}$ for $i = 1, \ldots, \alpha$ and then rename $a$ by $a_1$ we get
$a_i \notin \acl^n(\{a_1, \ldots, a_{i-1}\} \cup B)$ for all $i  = 1, \ldots, \alpha + 1$.

Now suppose that $a_1, \ldots, a_{\alpha+1} \in \acl^n(A)$ and, for all $k = 1, \ldots, \alpha + 1$,
\[
a_k \notin \acl^n(\{a_1, \ldots, a_{k-1}\} \cup B).
\]
By the induction hypothesis we have 
$\rk^n(A / a_1 B) \geq \alpha$.
By the choice of $a_1, \ldots, a_{\alpha + 1}$ we have $a_1 \in \acl^n(A) \setminus \acl^n(B)$,
so $\rk^n(A / B) \geq \alpha + 1$.
Now we have proved part~(i).
Part~(ii) follows directly from part~(i).

(iii) Suppose that $\rk^n(A / B) = \alpha$ and that
$a_1, \ldots, a_\alpha \in \acl^n(A)$ are such that, for all $k = 1, \ldots, \alpha$,
$a_k \notin \acl^n(\{a_1, \ldots, a_{k-1}\} \cup B)$. 
Then $\rk^n(a_k / \{a_1, \ldots, a_{k-1}\} \cup B) \geq 1$ for all $k = 1, \ldots, \alpha$.
Suppose, for a contradiction, that for some $k$, $\rk^n(a_k / \{a_1, \ldots, a_{k-1}\} \cup B) \geq 2$.
Choose the least such $k$.
Then there is 
\[
a \in \acl^n(a_k) \setminus \acl^n(\{a_1, \ldots, a_{k-1}\} \cup B)
\]
such that $\rk^n(a_k / \{a, a_1, \ldots, a_{k-1}\} \cup B) \geq 1$
which in particular means that 
\[
a_k \notin \acl^n(\{a, a_1, \ldots, a_{k-1}\} \cup B).
\]
Since $a \in \acl^n(a_k)$ it follows from the choice of $a_1, \ldots, a_\alpha$ that, 
for all $i = k+1, \ldots, \alpha$,
$a_i \notin \acl^n(\{a, a_1, \ldots, a_{i-1}\} \cup B)$.
Since $a \in \acl^n(a_k) \subseteq \acl^n(A)$ it now follows from part~(i) that $\rk^n(A / B) \geq \alpha + 1$,
contradicting the assumption.

By assumption and part~(ii), for every $k = 1, \ldots, \alpha$, 
the number $\alpha - k$ is maximal such that there are $a'_{k+1}, \ldots, a'_\alpha \in \acl^n(A)$
such that, for all $i = k+1, \ldots, \alpha$, 
\[
a'_i \notin \acl\meq(\{a_1, \ldots, a_k, a'_{k+1}, \ldots, a'_{i-1}\} \cup B).
\]
Hence $\rk(A / \{a_1, \ldots, a_k\} \cup B) = \alpha - k$ for all $k = 1, \ldots, \alpha$.

Let $k \in \{1, \ldots, \alpha\}$.
By part~(i), the sequence $a_1, \ldots, a_k$ witnesses that $\rk(a_1, \ldots, a_k / B)$ $\geq k$.
Suppose, for a contradiction, that $\rk(a_1, \ldots, a_k / B) \geq k+1$.
Then (by part ~(i)) we find 
\[
a'_0, \ldots, a'_k \in \acl^n(a_1, \ldots, a_k) \subseteq \acl^n(A)
\]
such that, for all $i = 0, \ldots, k$, $a'_i \notin \acl^n(\{a'_1, \ldots, a'_{i-1}\} \cup B)$.
Then $\acl^n(\{a'_0, \ldots, a'_k\} \cup B) \subseteq \acl^n(\{a_1, \ldots, a_k\} \cup B)$
so, for all $i = k+1, \ldots, \alpha$, 
$a_i \notin \acl^n(\{a'_0, \ldots, a'_k, a_{k+1}, \ldots, a_{i-1}\} \cup B)$
By part~(i), the sequence $a'_0, \ldots, a'_k, a_{k+1}, \ldots, a_\alpha$ witnesses that $\rk(A / B) \geq \alpha+1$ 
which contradicts the assumption.
\hfill $\square$

\begin{defin}\label{definition of n-cs}{\rm
Let $A, B \subseteq M\meq$ and suppose that $\rk^n(A / B) = \alpha < \omega$.
Then every sequence
$a_1, \ldots, a_\alpha \in \acl^n(A)$ such that, for all $k = 1, \ldots, \alpha$,
$a_k \notin \acl^n(\{a_1, \ldots, a_{k-1}\} \cup B)$
will be called an {\em $n$-coordinatization sequence} ($n$-cs) for $A/B$ (``$A$ over $B$'').
}\end{defin}

\begin{lem}\label{ranks of finite sets are finite}
If $A, B \subseteq M\meq$ and $A$ is finite then $\rk^n(A / B)$ is defined and finite.
\end{lem}

\noindent
{\bf Proof.}
Let $A, B \subseteq M\meq$ where $A$ is finite. 
Then $|\acl^n(A)| = \alpha$ for some $\alpha < \omega$.
If $\rk^n(A / B) \geq \alpha + 1$ then, by
Lemma~\ref{characterization of rank}~(i),
there are $a_1, \ldots, a_{\alpha + 1} \in \acl^n(A)$ such that, for all $k = 1, \ldots, \alpha + 1$,
$a_k \notin \acl^n(\{a_1, \ldots, a_{k-1}\} \cup B)$.
Then $a_i \neq a_j$ if $i \neq j$ so $|\acl^n(A)| \geq \alpha + 1$, a contradiction.
\hfill $\square$

\begin{lem}\label{reduction of n-rank to 0-rank}
Suppose that $T$ has geometric elimination of imaginaries.
Let $\alpha < \omega$ and $A, B \subseteq M\meq$.\\
(i) If $\rk^n(A / B) \geq \alpha$ then $\rk^0(A / B) \geq \alpha$.\\
(ii) If $A$ is finite then $\rk^n(A / B) = \rk^0(A / B)$.
\end{lem}

\noindent
{\bf Proof.}
(i) Suppose that $\rk^n(A / B) \geq \alpha$.
By Lemma~\ref{characterization of rank},
there are $a_1, \ldots, a_\alpha \in \acl^n(A)$ such that, for all $k = 1, \ldots, \alpha$,
$a_k \notin \acl^n(\{a_1, \ldots, a_{k-1}\} \cup B)$.
By induction on $\alpha$ we prove that there are $a'_1, \ldots, a'_\alpha \in \acl^0(A)$ 
such that $a'_k \notin \acl^0(\{a'_1, \ldots, a'_{k-1}\} \cup B)$ 
and $a'_k \in \acl^0(a_k)$ for all $k = 1, \ldots, \alpha$,
and it follows from 
Lemma~\ref{characterization of rank} 
that $\rk^0(A / B) \geq \alpha$.
If $\alpha = 0$ there is nothing to prove, so suppose that $\alpha > 0$.
By the induction hypothesis there are $a'_1, \ldots, a'_{\alpha-1} \in \acl^0(A)$ such that
$a'_k \notin \acl^0(\{a'_1, \ldots, a'_{k-1}\} \cup B)$ 
and $a'_k \in \acl^0(a_k)$ for all $k = 1, \ldots, \alpha-1$.
Since $T$ is assumed to have geometric elimination of imaginaries, there is
a finite sequence $\bar{c}$ of elements from $M_0 = M$ such that $\acl\meq(a_\alpha) = \acl\meq(\bar{c})$.
Since $a_\alpha \notin \acl^n(\{a_1, \ldots, a_{\alpha-1}\} \cup B)$
we have $\rng(\bar{c}) \not\subseteq  \acl^n(\{a_1, \ldots, a_{\alpha-1}\} \cup B)$.
As $a'_k \in \acl^0(a_k)$ for all $k = 1, \ldots, \alpha - 1$ it follows that
$\rng(\bar{c}) \not\subseteq \acl^0(\{a'_1, \ldots, a'_{\alpha - 1}\} \cup B)$.
Thus we can choose some element from $\bar{c}$, which we call $a'_\alpha$, 
such that $a'_\alpha \notin \acl^0(\{a'_1, \ldots, a'_{\alpha-1}\} \cup B)$.
Now Lemma~\ref{characterization of rank} implies that $\rk^0(A / B) \geq \alpha$.

(ii) Suppose that $A$ is finite.
By Lemma~\ref{ranks of finite sets are finite}, 
$\rk^n(A / B) = \alpha$ for some $\alpha < \omega$.
By Lemma~\ref{n+1-rank is at least as big as n-rank},
$\rk^n(A / B) \geq \rk^0(A / B)$,
and by part~(i)
$\rk^n(A / B) \leq \rk^0(A / B)$.
\hfill $\square$

\begin{lem}\label{monotonicity of n-rank}
Let $A, B, C, D \subseteq M\meq$ where $A \supseteq B$, $C \subseteq D$ and $A$ is finite.
Then
\begin{itemize}
\item[(a)] $\rk^n(A / A) = 0$,
\item[(b)] $\rk^n(A / C) \geq \rk^n(B / C)$,
\item[(c)] $\rk^n(A / C) \geq \rk^n(A / D)$, and
\item[(d)] $\rk^n(A / C) \geq \rk^n(B / C) + \rk^n(A / BC)$.
\end{itemize}
\end{lem}

\noindent
{\bf Proof.}
Suppose that $A, B, C, D \subseteq M\meq$ where $A \supseteq B$, $C \subseteq D$ and $A$ is finite,
so all mentioned $n$-ranks in (a) -- (d) are finite.
Part~(a) is obvious from the definition of $n$-rank.

For part~(b), suppose that $\rk^n(B / C) \geq \alpha$ (where $\alpha < \omega$).
By Lemma~\ref{characterization of rank},
there are $a_1, \ldots, a_\alpha \in \acl^n(B)$ such that, for all $k = 1, \ldots, \alpha$,
$a_k \notin \acl^n(\{a_1, \ldots, a_{k-1}\} \cup C)$.
Since $A \supseteq B$ it follows that $a_1, \ldots, a_\alpha \in \acl^n(A)$, so, by
Lemma~\ref{characterization of rank} again,
$\rk^n(A / C) \geq \alpha$.

Now consider part~(c). Suppose that $\rk^n(A / D) \geq \alpha$.
Then there are $a_1, \ldots, a_\alpha \in \acl^n(A)$ such that, for all $k = 1, \ldots, \alpha$,
$a_k \notin \acl^n(\{a_1, \ldots, a_{k-1}\} \cup D)$.
Since $C \subseteq D$ we get $a_k \notin \acl^n(\{a_1, \ldots, a_{k-1}\} \cup C)$ for all $k = 1, \ldots, \alpha$,
so $\rk^n(A / C) \geq \alpha$.

Finally we consider part (d). 
Let $\beta = \rk^n(B / C)$ and $\alpha = \rk^n(A / BC)$, so $\alpha, \beta < \omega$.
By Lemma~\ref{characterization of rank}, there are $b_1, \ldots b_\beta \in \acl^n(B)$ and
$a_1, \ldots, a_\alpha \in \acl^n(A)$ such that, for all $k = 1, \ldots, \beta$,
$b_k \notin \acl^n(\{b_1, \ldots, b_{k-1}\} \cup C)$, and for all $l = 1, \ldots, \alpha$,
$a_l \notin \acl^n(\{a_1, \ldots, a_{l-1}\} \cup BC)$.
Moreover, $\beta$ and $\alpha$ are maximal such that such sequences exist.
For a contradiction, suppose that, for some $l \in \{1, \ldots, \alpha\}$,
$a_l \in \acl^n(\{b_1, \ldots, b_\beta, a_1, \ldots, a_{l-1}\} \cup C)$.
Since $b_1, \ldots, b_\beta \in \acl^n(B)$ we get
$a_l \in \acl^n(\{a_1, \ldots, a_{l-1}\} \cup BC)$ which contradicts the choice of $a_1, \ldots, a_\alpha$.
Hence we conclude that $a_l \notin \acl^n(\{b_1, \ldots, b_\beta, a_1, \ldots, a_{l-1}\} \cup C)$
for all $l = 1, \ldots, \alpha$.
Since $B \subseteq A$ we have $\acl^n(B) \subseteq \acl^n(A)$ and therefore
$b_1, \ldots, b_\beta, a_1, \ldots, a_\alpha \in \acl^n(A)$.
By Lemma~\ref{characterization of rank},
we get $\rk(A / C) \geq \beta + \alpha$.
\hfill $\square$

\begin{lem}\label{finding a in acl of rank 1 over something else}
Let $A, B \subseteq M\meq$.
If $\acl^n(A) \setminus \acl^n(B) \neq \es$ then there is 
$a \in \acl^n(A) \setminus \acl^n(B)$ such that $\rk^n(a / B) = 1$.
\end{lem}

\noindent
{\bf Proof.}
Note that (by Lemma~\ref{ranks of finite sets are finite}),
if $a \in \acl^n(A) \setminus \acl^n(B)$ then $\rk^n(a / B) < \omega$.
We first prove the following claim:

\medskip

{\bf Claim.} If $a' \in \acl^n(A) \setminus \acl^n(B)$ and $\rk^n(a' / B) = \alpha + 1$ where $\alpha \geq 1$,
then there is $a \in \acl^n(a') \subseteq \acl^n(A)$ such that $1 \leq \rk^n(a / B) \leq \alpha$.

\medskip

\noindent
{\em Proof of the claim.}
Suppose that $a' \in \acl^n(A) \setminus \acl^n(B)$ and $\rk^n(a' / B) = \alpha + 1$ where $\alpha \geq 1$.
Suppose that $a' \in \acl^n(A) \setminus \acl^n(B)$ and $\rk^n(a' / B) = \alpha + 1$
where $\alpha \geq 1$.
By the definition of $\rk^n$ there is $a \in \acl^n(a') \setminus \acl^n(B)$ such that $\rk^n(a' / aB) = \alpha$.
If $a' \in \acl^n(aB)$ then $\rk^n(a' / aB) = 0$, so $\alpha = 0$, which contradicts that $\alpha \geq 1$.
Hence $a' \notin \acl^n(aB)$.
As $a \notin \acl^n(B)$ we have $\rk^n(a / B) \geq 1$.

We now show that $\rk^n(a / B) \leq \alpha$.
For a contradiction, suppose that $\rk^n(a / B) \geq \alpha + 1$.
By Lemma~\ref{characterization of rank} 
there are $a_1, \ldots, a_{\alpha + 1} \in \acl^n(a)$ such that for all $i = 1, \ldots, \alpha + 1$,
$a_i \notin \acl^n(\{a_1, \ldots, a_{i-1}\}  \cup B)$.
If $a' \in \acl^n(\{a_1, \ldots, a_{\alpha + 1}\} \cup B)$,
then $a' \in \acl^n(aB)$ contradicting what we concluded above.
Hence $a' \notin \acl^n(aB)$.
Therefore the sequence $a_1, \ldots, a_{\alpha + 1}, a' \in \acl^n(a')$ and 
Lemma~\ref{characterization of rank}
witness that $\rk^n(a' / B) \geq \alpha + 2$ which contradicts the assumption about $a'$.
This concludes the proof of the claim.

Now suppose that $a' \in \acl^n(A) \setminus \acl^n(B)$.
Let $\alpha = \rk^n(a' / B)$. Since $a' \notin \acl^n(B)$ we have $\alpha \geq 1$.
If $\alpha = 1$ then we are done (by letting $a = a'$).
So suppose that $\alpha \geq 2$.
By repeatedly using the claim we eventually find $a \in \acl^n(a') \subseteq \acl^n(A)$ such that $\rk^n(a / B) = 1$.
This ends the proof of the lemma.
\hfill $\square$

\begin{lem}\label{locality property of rank}
(i) Suppose that $A \subseteq M\meq$ is finite and $B \subseteq C \subseteq M\meq$.
If $\rk^n(A / C) < \rk^n(A / B)$ then there is a finite $C' \subseteq C$ such that $\rk^n(A / B \cup C') < \rk^n(A / B)$.\\
(ii) Suppose that $A \subseteq M\meq$ is finite and $B \subseteq M\meq$.
Then there is finite $B' \subseteq B$ such that $\rk^n(A / B') = \rk^n(A / B)$.
\end{lem}

\noindent
{\bf Proof.}
(i) Let $A \subseteq M\meq$ be finite and let $B \subseteq C \subseteq M\meq$.
We prove the claim by induction on $\rk^n(A / B)$.
If $\rk^n(A / B) = 0$ then there is nothing to prove (as we cannot have $\rk^n(A / C) < \rk^n(A / B)$).
So suppose that $\rk^n(A / B) = \alpha + 1$.
Then there is $a \in \acl^n(A) \setminus \acl^n(B)$ such that $\rk^n(A / \{a\} \cup B) \geq \alpha$.
In fact we must have $\rk^n(A / \{a\} \cup B) = \alpha$ because if $\rk^n(A / \{a\} \cup B) \geq \alpha+1$ then $
\rk^n(A / B) \geq \alpha+2$ which contradicts the assumption.
Since $A$ is finite it follows that  $\acl^n(A)$ is finite and therefore there is finite $C_1 \subseteq C$ such that
\[
\acl^n(A) \cap \acl^n(C) = \acl^n(A) \cap \acl^n(C_1).
\]
By the induction hypothesis, if $a \in \acl^n(A) \setminus \acl^n(B)$ is such that 
$\rk^n(A / \{a\} \cup C) < \rk^n(A / \{a\} \cup B) = \alpha$,
then there is finite $C_{a} \subseteq C$ such that $\rk^n(A / \{a\} \cup B \cup C_{a}) < \rk^n(A / \{a\} \cup B)$.
Let $C_2$ be the union of all $C_{a} \subseteq C$ where $a$ ranges over the members of the finite set $\acl^n(A)$
such that $\rk^n(A / \{a\} \cup C) < \rk^n(A / \{a\} \cup B)$.
Let $C' = C_1 \cup C_2$ so $C'$ is finite.

Now suppose that $\rk^n(A / C) < \rk^n(A / B)$. Then {\em for every} $a \in \acl^n(A) \setminus \acl^n(B)$ such that
 $\rk^n(A / \{a\} \cup B) = \alpha$ we have 
either $a \in \acl\meq(C)$ or $\rk^n(A / \{a\} \cup C) < \alpha$ which in turn implies that either
$a \in \acl^n(C')$ or $\rk^n(A / \{a\} \cup B \cup C') < \rk^n(A / \{a\} \cup B)$.
It follows that $\rk^n(A / B \cup C') < \alpha+1 = \rk^n(A / B)$.

(ii) Suppose that $\rk^n(A / B) = \alpha$ (where $A$ is finite) and $\rk^n(A / \es) = \beta$ so $\alpha \leq \beta$
(by Lemma~\ref{monotonicity of n-rank}).
If $\alpha = \beta$ we are done. So suppose that $\alpha < \beta$. Let $B_0 = \es$.
By part~(i) there is finite $B_1 \subseteq B$ such that $\rk^n(A / B_1) < \rk^n(A / B_0)$.
If $\rk^n(A / B) < \rk^n(A / B_1)$ then we use part~(i) again and get a finite $B_2 \subseteq B$ such that
$\rk^n(A / B_0 \cup B_1  \cup B_2) < \rk^n(A / B_0  \cup B_1)$. We can repeat this procedure as long as 
$\rk^n(A / B) < \rk^n(A / B_0 \cup \ldots \cup B_m)$. But as there is no infinite decreasing sequence of natural numbers
we eventually find $m <  \omega$ and finite $B_k \subseteq B$ 
for $k \leq m$ such that $\rk^n(A / B_0 \cup \ldots \cup B_m) = \rk^n(A / B)$.
\hfill $\square$

\begin{lem}\label{locality property again}
Suppose that $A \subseteq M\meq$ is finite and $B \subseteq M\meq$.
Then there is countable $C \subseteq B$ such that $\rk^n(A / C) = \rk^n(A / B)$ for all $n < \omega$.
\end{lem}

\noindent
{\bf Proof.}
Suppose that $A \subseteq M\meq$ is finite and $B \subseteq M\meq$.
Lemma~\ref{locality property of rank} says that for every $n < \omega$ there is a finite $B_n \subseteq B$
such that $\rk^n(A / B_n) = \rk^n(A / B)$. 
Let $C = \bigcup_{n < \omega} B_n$, so $C \subseteq B$ is countable.
By Lemma~\ref{monotonicity of n-rank}, for all $n < \omega$,
$\rk^n(A / B_n) \geq \rk^n(A / C) \geq \rk^n(A / B) = \rk^n(A / B_n)$, and hence
$\rk^n(A / C) = \rk^n(A / B_n) = \rk^n(A / B)$.
\hfill $\square$

\begin{lem}\label{n-rank is determined by type}
(i) Let $m, r, k, l < \omega$ where $n \leq m$. The $(k+l)$-ary relation on $M_m$ which holds for
$(a_1, \ldots, a_k, b_1, \ldots, b_l) \in (M_m)^{k+l}$ if and only if 
$\rk^n(a_1, \ldots, a_k / b_1, \ldots, b_l) = r$ is $\es$-definable in $\mcM_m$ and in $\mcM\meq$.\\
(ii) Let $\bar{a}$, $\bar{a}'$, $\bar{b}$ and $\bar{b}'$ be finite sequences of elements from $M\meq$.
If $\tp(\bar{a}, \bar{b}) = \tp(\bar{a}', \bar{b}')$, then $\rk^n(\bar{a} / \bar{b}) = \rk^n(\bar{a}' / \bar{b}')$.
\end{lem}

\noindent
{\bf Proof.}
(i) Suppose that $n \leq m < \omega$.
Since $T$ is $\omega$-categorical it follows that $Th(\mcM_m)$ is $\omega$-categorical and from this it follows that,
for all $k < \omega$, the $(k+1)$-ary relation on $M_m$ which holds
for $(b, a_1, \ldots, a_k) \in (M_m)^{k+1}$ if and only if $b \in \acl^n(a_1, \ldots, a_k)$ 
is $\es$-definable in $\mcM_m$.
It follows that for all $r, k, l < \omega$, the $(k+l)$-ary relation on $M_m$ which holds for 
$(a_1, \ldots, a_k, b_1, \ldots, b_l) \in (M_m)^{k+l}$
if and only if $r$ is maximal  such that 
\begin{align}\label{the condition for rank being r}
&\text{there are  } \ a'_1, \ldots, a'_r \in \acl^n(a_1, \ldots, a_k) \ \text{ such that,}\\
&\text{for all $i = 1, \ldots, r$,} \
a'_i \notin \acl^n(a'_1, \ldots, a'_{i-1}, b_1, \ldots, b_l) \nonumber
\end{align}
is $\es$-definable in $\mcM_m$, by $\varphi_r^m(x_1, \ldots, x_k, y_1, \ldots, y_l)$ say.
Then the formula $\psi_r^m(x_1, \ldots, x_k,$ $y_1, \ldots, y_l)$ which expresses that 
``all $x_1, \ldots, x_k, y_1, \ldots, y_l$ belong to $M_m$ and $\varphi_r^m(x_1, \ldots, x_k,$ $y_1, \ldots, y_l)$ holds''
defines the same relation in $\mcM\meq$.

By 
Lemma~\ref{characterization of rank},
$\rk^n(a_1, \ldots, a_k / b_1, \ldots, b_l) = r$ if and only if $r$ is maximal such
that~(\ref{the condition for rank being r}) holds. It follows that, for every $r < \omega$,
the $(k+l)$-ary relation on $(M_m)^{k+l}$ which holds for 
$(a_1, \ldots, a_k, b_1, \ldots, b_l) \in (M_m)^{k+l}$
if and only if  $\rk^n(a_1, \ldots, a_k / b_1, \ldots, b_l) = r$
is $\es$-definable in $\mcM_m$ by $\varphi_r^m(x_1, \ldots, x_k, y_1, \ldots, y_l)$.
The same relation is definable in $\mcM\meq$ by $\psi_r^m(x_1, \ldots, x_k, y_1, \ldots, y_l)$.

(ii) Suppose that $\bar{a}, \bar{a}' \in (M\meq)^k$, $\bar{b}, \bar{b}' \in (M\meq)^l$ and that
$\rk^n(\bar{a} / \bar{b}) = r$. 
Then there is $m < \omega$ such that $n \leq m$, $\bar{a}, \bar{a}' \in (M_m)^k$, 
and $\bar{b}, \bar{b}' \in (M_m)^l$.
Then $\mcM\meq \models \psi_r^m(\bar{a}, \bar{b})$, and if $\tp(\bar{a}, \bar{b}) = \tp(\bar{a}', \bar{b}')$
then also $\mcM\meq \models \psi_r^m(\bar{a}', \bar{b}')$, so $\rk^n(\bar{a}' / \bar{b}') = r$.
\hfill $\square$

\begin{lem}\label{nonalgebraic extension}
Let $B \subseteq C \subseteq M\meq$ be finite and $a \in M\meq$.
If $a \notin \acl\meq(B)$ then there is $a' \in M\meq$ such that $\tp(a' / B) = \tp(a / B)$
and $a' \notin \acl\meq(C)$.
\end{lem}

\noindent
{\bf Proof.}
Let $S$ be the sort of $a$.
Since $a \notin \acl\meq(B)$ it follows that $\tp(a / B)$ has infinitely many realizations.
As $C$ is finite it follows (from $\omega$-categoricity) that $\acl\meq(C) \cap S$ is finite and hence there is $a' \notin \acl\meq(C)$
such that $\tp(a' / B) = \tp(a / B)$.
\hfill $\square$

\begin{lem}\label{n-rank extension}
Suppose that $\bar{d}$ is a finite sequence of elements from $M_n$, 
$B, C \subseteq M\meq$ are finite and $B \subseteq C$.
Let $\alpha = \rk^n(\bar{d} / B)$ (so $\alpha < \omega$).
Then there is a finite sequence $\bar{d}'$ of elements from $M_n$ 
such that $\tp(\bar{d}' / B) = \tp(\bar{d} / B)$ and $\rk^n(\bar{d}' / C) = \alpha$.
\end{lem}

\noindent
{\bf Proof.}
Let $\bar{d}$ be a finite sequence of elements from $M_n$.
Suppose that $\alpha = \rk^n(\bar{d} / B)$.
By Lemma~\ref{characterization of rank}
there are $a_1, \ldots, a_\alpha \in \acl^n(\bar{d})$ such that,
\begin{align}\label{characterization of rank in proving n-extension}
\text{ for all $k = 1, \ldots, \alpha$, } \ &\rk^n(a_k / \{a_1, \ldots, a_{k-1}\} \cup B) = 1, \\
&\rk^n(a_1, \ldots, a_k / B) = k, \text{ and} \nonumber \\
&\rk^n(\bar{d} / \{a_1, \ldots, a_k\} \cup B) = \alpha-k. \nonumber
\end{align}
In particular, 
$\rk^n(a_1, \ldots, a_\alpha / B) = \alpha$ and $\rk^n(\bar{d} / \{a_1, \ldots, a_\alpha\} \cup B) = 0$,
so (as $\rng(\bar{d}) \subseteq M_n$) 
$\rng(\bar{d}) \subseteq \acl^n(\bar{d}) \subseteq \acl^n(\{a_1, \ldots, a_\alpha\} \cup B)$ and hence
$\acl^n(\bar{d}B) = \acl^n(\{a_1, \ldots, a_\alpha\} \cup B)$.

Suppose that $k < \alpha$ and that there are $a'_1, \ldots, a'_k \in M_n$ such that 
\begin{align}\label{a' have same type as a}
&\tp(a'_1, \ldots, a'_k / B) = \tp(a_1, \ldots, a_k / B), \  \text{ and} \\
\label{a'-i+1 not in acl of}
&\text{for all } i = 1, \ldots, k,  \ 
a'_i \notin \acl^n(\{a'_1, \ldots, a'_{i-1}\} \cup C).
\end{align}
We will find $a'_{k+1}$ so that the above holds also with $k$ replaced by $k+1$.
Let 
\[
p(x, y_1, \ldots, y_k) = \tp(a_{k+1}, a_1, \ldots, a_k / B).
\]
By ~(\ref{characterization of rank in proving n-extension}),
$p(x, a_1, \ldots, a_k)$ is non-algebraic.
By~(\ref{a' have same type as a}), also $p(x, a'_1, \ldots, a'_k)$ is non-algebraic, so
by Lemma~\ref{nonalgebraic extension}
there is $a'_{k+1} \in M_n \setminus \acl^n(\{a'_1, \ldots, a'_k\} \cup C)$  which realizes $p(x, a'_1, \ldots, a'_k)$, 
that is, $\tp(a'_1, \ldots, a'_k, a'_{k+1} / B) = \tp(a_1, \ldots, a_k, a_{k+1} / B)$.
Now~(\ref{a' have same type as a}) and~(\ref{a'-i+1 not in acl of}) hold if $k$ is replaced by $k+1$.
By induction it follows that there are $a'_1, \ldots, a'_\alpha \in M_n$ such that
\begin{align}\label{a' have same type as a, for alpha}
&\tp(a'_1, \ldots, a'_\alpha / B) = \tp(a_1, \ldots, a_\alpha / B), \  \text{ and} \\
\label{a'-i+1 not in acl of, for alpha}
&\text{for all } i = 1, \ldots, \alpha,  \ 
a'_i \notin \acl^n(\{a'_1, \ldots, a'_{i-1}\} \cup C).
\end{align}
From~(\ref{a' have same type as a, for alpha}) it follows that there is $\bar{d}'$ such that
\begin{equation}\label{A and A' have the same type}
\tp(\bar{d}', a'_1, \ldots, a'_\alpha / B) = \tp(\bar{d}, a_1, \ldots, a_\alpha / B).
\end{equation}
Then $a'_1, \ldots, a'_\alpha \in \acl^n(\bar{d}')$.
This together with~(\ref{a'-i+1 not in acl of, for alpha}) implies that
$\rk^n(\bar{d}' / C) \geq \alpha$.
From~(\ref{A and A' have the same type})
and Lemma~\ref{n-rank is determined by type}
we get $\rk^n(\bar{d}' / B) = \rk^n(\bar{d} / B) = \alpha$.
Since we must
(by Lemma~\ref{monotonicity of n-rank}) 
have $\rk^n(\bar{d}' / C) \leq \rk^n(\bar{d}' /B)$ we get $\rk^n(\bar{d}' / C) = \alpha$.
From~(\ref{A and A' have the same type}) we get $\tp(\bar{d}' / B) = \tp(\bar{d} / B)$ so the proof is completed.
\hfill $\square$

\begin{lem}\label{general n-rank extension}
Suppose that $\bar{d}$ is a finite sequence of elements from $M_n$, $B, C \subseteq M\meq$ and $B \subseteq C$.
Then there is a finite sequence $\bar{d}'$ of elements from $M_n$ such that 
$\tp(\bar{d}' / B) = \tp(\bar{d} / B)$ and $\rk^n(\bar{d}' / C) = \rk^n(\bar{d} / B)$.
\end{lem}

\noindent
{\bf Proof.}
Let $\bar{d}$ be a finite sequence of elements from $M_n$ and let $\alpha = \rk^n(\bar{d} / B)$.
Hence $\rk^n(\bar{d} / B') \geq \alpha$ for every finite $B' \subseteq B$.
By Lemma~\ref{n-rank is determined by type},
for all finite $B' \subseteq B$ and finite $C' \subseteq C$
there is a formula $\varphi_{B'C'}(\bar{x})$ with parameters from $B' \cup C'$ that expresses that
$\tp(\bar{x} / B') = \tp(\bar{d} / B')$ and
$\rk^n(\bar{x} / B'C') \geq \alpha$.
Let $\Phi(\bar{x})$ be the set of all such formulas $\varphi_{B'C'}(\bar{x})$ as $B'$ and $C'$
varies over finite subsets of $B$ and $C$, respectively.
From Lemma~\ref{n-rank extension}
it follows that every finite subset of $\Phi(\bar{x})$ is consistent.
From compactness it follows that $\Phi(\bar{x})$ is consistent.
Let $\bar{d}'$ be a realization of the type $\Phi(\bar{x})$.
It follows that $\tp(\bar{d}' / B) = \tp(\bar{d} / B)$ and, by Lemma~\ref{locality property of rank}, that
$\rk^n(\bar{d}' / C) \geq \alpha$.
Since we must 
(by Lemmas~\ref{monotonicity of n-rank} and~\ref{n-rank is determined by type}) 
have $\rk^n(\bar{d}' / C) \leq \rk^n(\bar{d}' / B) = \rk^n(\bar{d} / B) = \alpha$
it follows that $\rk^n(\bar{d}' / C) = \alpha$.
\hfill $\square$

\section{$n$-independence}\label{n-independence}

\noindent
Now we are ready to define a hierarchy of independence relations, the notions of $n$-independence, 
for all $n < \omega$.
In this section we will prove that, for all $n < \omega$, $n$-independence has all properties
of an independence relation restricted to $M_n$ with the {\em possible exception} of the symmetry property.
Actually the restriction to $M_n$ is only needed for the extension property.
We also show that if the theory $T$ has geometric elimination of imaginaries, then the hierarchy of
$n$-independence relations collapses to the bottom level of 0-independence.

\begin{defin}\label{definition of independence}{\rm
Let $n < \omega$.
For $A, B, C \subseteq M\meq$ define
$A \underset{C}{\indn} B$ if and only if, for all finite $A' \subseteq A$, $\rk^n(A' / BC) = \rk^n(A' / C)$.
If $\bar{a}$ and $\bar{b}$ are sequences then $\bar{a} \underset{C}{\indn} \bar{b}$ means the same as 
$\rng(\bar{a}) \underset{C}{\indn} \rng(\bar{b})$.
}\end{defin}

\begin{exam}\label{example for n-independence}{\rm
Let $T$ be as in Example~\ref{example for n-rank}, $\mcM \models T$ and let $a, b \in M$ be distinct
and such that $[a]_E = [b]_E$. By the conclusions in that example we get
$\rk^0(a / b) = \rk^0(a)$ and $\rk^1(a / b) = 1 < 2 = \rk^1(a)$.
Hence $a \underset{\es}{\indo} b$ and $a \underset{\es}{\nindone} b$.
}\end{exam}

\noindent
{\em For the rest of this section we fix an arbitrary $n < \omega$.}

\begin{lem}\label{invariance} {\bf (Invariance)} 
Suppose that $\sigma$ is an elementary function from a subset of $M\meq$ to a subset of $M\meq$
and that $A, B$, and $C$
are subsets of the domain of $\sigma$.
If $A \underset{C}{\indn} B$ then $\sigma(A) \underset{\sigma(C)}{\indn} \sigma(B)$.
\end{lem}

\noindent
{\bf Proof.}
Suppose that $\sigma$ is an elementary function and that $A, B$, and $C$
are subsets of the domain of $\sigma$.
Suppose that $A \underset{C}{\indn} B$, so
by Definition~\ref{definition of independence} of $\indn$
we have $\rk^n(A' / BC) = \rk^n(A' / C)$ for all finite $A' \subseteq A$.
It suffices to show that for every finite $A' \subseteq A$, 
$\rk^n(\sigma(A') / \sigma(B)\sigma(C)) = \rk^n(\sigma(A') / \sigma(C))$.
This follows if we can show that for all $A, B \subseteq M\meq$ where $A$ is finite,
$\rk^n(A / B) = \rk^n(\sigma(A) / \sigma(B))$.
Since also the inverse of $\sigma$ is an elementary function, it actually suffices to show that
$\rk^n(A / B) \geq \rk^n(\sigma(A) / \sigma(B))$.

So let $A = \{a_1, \ldots, a_k\} \subseteq A$ be finite.
Lemma~\ref{locality property of rank}
tells that there is finite $B' = \{b_1, \ldots, b_l\} \subseteq B$ such that 
$\rk^n(A / B') = \rk^n(A / B)$.
Since $\tp(a_1, \ldots, a_k, b_1, \ldots, b_l) = \tp(\sigma(a_1), \ldots,  \sigma(a_k), \sigma(b_1), \ldots, \sigma(b_l))$
it follows from 
Lemmas~\ref{monotonicity of n-rank} and~\ref{n-rank is determined by type}
that 
\[
\rk^n(\sigma(A) / \sigma(B)) \leq \rk^n(\sigma(a_1), \ldots,  \sigma(a_k)/ \sigma(b_1), \ldots, \sigma(b_l)) = 
\rk^n(A / B') = \rk^n(A / B).
\]
\hfill $\square$

\begin{lem}\label{monotonicity of independence} {\bf (Monotonicity)}
Let $A, B, C, D \subseteq M\meq$ where $B \subseteq C \subseteq D$.\\
If $A \underset{B}{\indn} D$, then $A \underset{B}{\indn} C$ and $A \underset{C}{\indn} D$.
\end{lem}

\noindent
{\bf Proof.}
Let $A, B, C, D \subseteq M\meq$ where $B \subseteq C \subseteq D$.
Suppose that $A \underset{B}{\indn} D$ and let $A' \subseteq A$ be finite.
Then $\rk^n(A' / D) = \rk^n(A' / B)$.
By Lemma~\ref{monotonicity of n-rank},
$\rk^n(A' / B) = \rk^n(A' / D) \leq \rk^n(A' / C) \leq \rk^n(A' / B)$.
Hence $\rk^n(A' / D) = \rk^n(A' / C) = \rk^n(A' / B)$ and thus $A \underset{B}{\indn} C$
and $A \underset{C}{\indn} D$.
\hfill $\square$

\begin{lem}\label{transitivity of independence} {\bf (Transitivity)}
Suppose that $A, B, C, D \subseteq M\meq$ and $B \subseteq C \subseteq D$.
If $A \underset{B}{\indn} C$ and $A \underset{C}{\indn} D$ then $A \underset{B}{\indn} D$.
\end{lem}

\noindent
{\bf Proof.}
Suppose that $A, B, C, D \subseteq M\meq$ and $B \subseteq C \subseteq D$.
By the definition of $n$-independence it suffices to prove that for all {\em finite} $A' \subseteq A$,
if $A' \underset{B}{\indn} C$ and $A' \underset{C}{\indn} D$ then $A' \underset{B}{\indn} D$.
Therefore we can assume that $A$ is finite.
So suppose that $A \underset{B}{\indn} C$ and $A \underset{C}{\indn} D$.
Then $\rk^n(A / B) = \rk^n(A / C) = \rk^n(A / D)$, so 
$A \underset{B}{\indn} D$.
\hfill $\square$

\begin{lem}\label{finite character of independence} {\bf (Finite character)} 
Let $A, B, C \subseteq M\meq$ and suppose that $A \underset{C}{\nindn} B$.
Then there are finite $A' \subseteq A$ and finite $B' \subseteq B$ such that $A' \underset{C}{\nindn} B'$.
\end{lem}

\noindent
{\bf Proof.}
If $A \underset{C}{\nindn} B$ then (by definition of $\nindn$)
$A'  \underset{C}{\nindn} B$ for some finite $A' \subseteq A$.
Then $\rk^n(A' / BC) < \rk^n(A' / C)$.
By Lemma~\ref{locality property of rank}~(i) there is finite $B' \subseteq B$ such that
$\rk^n(A' / B'C)$ $<$ $\rk^n(A' / C)$.
Hence $A'  \underset{C}{\nindn} B'$.
\hfill $\square$

\begin{lem}\label{locality of independence} {\bf (Locality)} 
If $A, B \subseteq M\meq$ and $A$ is {\em finite} then there is finite $C \subseteq B$ such that 
$A \underset{C}{\indn} B$.
\end{lem}

\noindent
{\bf Proof.}
Let $A, B \subseteq M\meq$ and suppose that $A$ is finite.
From Lemma~\ref{locality property of rank}~(ii) 
it follows that there is a finite $C \subseteq B$ such that 
$\rk^n(A / C) = \rk^n(A / B)$.
Hence $A \underset{C}{\indn} B$.
\hfill $\square$

\begin{lem}\label{extension property of n-independence} {\bf (Extension)}
Let $\bar{a}$ be a finite sequence of elements from $M_n$, 
and let $B, C \subseteq M\meq$ where $B \subseteq C$.
Then there is $\bar{a}' \subseteq M\meq$ such that $\tp(\bar{a}' / B) = \tp(\bar{a} / B)$ and
$\bar{a}' \underset{B}{\indn} C$.
\end{lem}

\noindent
{\bf Proof.}
Under the given assumptions, 
Lemma~\ref{general n-rank extension}
implies that there are is a finite sequence $\bar{a}'$ of elements from $M_n$ such that 
$\tp(\bar{a}' / B) = \tp(\bar{a} / B)$ and $\rk^n(\bar{a}' / C) = \rk^n(\bar{a} / B)$.
Hence $\rk^n(\bar{a}' / B) = \rk^n(\bar{a} / B)$, so $\bar{a}' \underset{B}{\indn} C$.
\hfill $\square$

\begin{lem}\label{n-independence and geometric elimination of imaginaries}
Suppose that $T$ has geometric elimination of imaginaries.
For all $n < \omega$ and all $A, B, C \subseteq M\meq$ we have
$A \underset{C}{\indn} B$ if and only if $A \underset{C}{\indo} B$.
\end{lem}

\noindent
{\bf Proof.}
By definition of $\indn$, $A \underset{C}{\indn} B$ if and only if
$\rk^n(A' / BC) = \rk^n(A' / C)$ for all finite $A' \subseteq A$.
Under the assumption that $T$ has geometric elimination of imaginaries, 
Lemma~\ref{reduction of n-rank to 0-rank}
implies that
$\rk^n(A' / BC) = \rk^0(A' / BC)$ and $\rk^n(A' / C) = \rk^0(A' / C)$ for all finite $A' \subseteq A$.
Hence $A \underset{C}{\indn} B$ if and only if $A \underset{C}{\indo} B$.
\hfill $\square$

\section{The exchange property and symmetry}\label{The exchange property and symmetry}

\noindent
According to the results in Section~\ref{n-independence},
for each $n < \omega$, $\indn$ has all properties of an independence relation restricted to $M_n$
with the {\em possible exception} of the symmetry property.
In this section we prove that, for any $n < \omega$, 
the symmetry of $\indn$ is a consequence of a more restricted symmetry property, namely
that the following assumption holds:

\begin{assump}\label{assumption about exchange property} 
{\bf (Exchange property with respect to $n$)} {\rm
In this section we fix an arbitrary $n < \omega$.
We assume that if $C \subseteq M\meq$, $2 \leq k < \omega$, $a_1, \ldots, a_k \in M_n$,
$\rk^n(a_i / C) = 1$ for all $i = 1, \ldots, k$, and
$a_k \in \acl\meq(\{a_1, \ldots, a_{k-1}\} \cup C) \setminus \acl\meq(\{a_2, \ldots, a_{k-1}\} \cup C)$,
then $a_1 \in \acl\meq(\{a_2, \ldots, a_k\} \cup C)$.
}\end{assump}

\begin{lem}\label{n-rank 1 elements form a pregeometry}
Let $C \subseteq M\meq$, $X = \{d \in M_n : \rk^n(d / C) = 1\}$,
and for every $A \subseteq X$, let $\cl(A) = \acl^n(AC) \cap X$.
Then $(X, \cl)$ is a pregeometry.
\end{lem}

\noindent
{\bf Proof.}
It is well-known that for all $A \subseteq M\meq$,
$A \subseteq \acl\meq(A)$, $\acl\meq(\acl\meq(A)) = \acl\meq(A)$,
and if $a \in \acl\meq(A)$ then $a \in \acl\meq(A')$ for some finite $A' \subseteq A$.
It follows straightforwardly that $\cl$ has properties (1) -- (3) 
in Definition~\ref{definition of pregeometry}, even without Assumption~\ref{assumption about exchange property}.
Suppose that $A \subseteq X$, $b, c \in X$, and $b \in \cl(A \cup \{c\}) \setminus \cl(A)$.
Then there is finite $A' \subseteq A$ such that $b \in \cl(A' \cup \{c\}) \setminus \cl(A)$.
By Assumption~\ref{assumption about exchange property},
we get $c \in \cl(A \cup \{b\})$.
Hence $(X, \cl)$ is a pregeometry.
\hfill $\square$

\medskip

\noindent
Due to Lemma~\ref{n-rank 1 elements form a pregeometry},
whenever $C \subseteq M\meq$, 
$X = \{d \in M_n : \rk^n(d / C) = 1\}$, and
$A \subseteq X$ it makes sense to talk about a {\em basis} of $A$ and the {\em dimension} of $A$
with respect to $(X, \cl)$ where $\cl$ is defined as in Lemma~\ref{n-rank 1 elements form a pregeometry}.

\begin{defin}\label{definition of n-ccs}{\rm
Let $A, B \subseteq M\meq$ and suppose that $A$ is finite.
An {\em $n$-canonical coordinatization sequence} ($n$-ccs) for $A/B$ (``$A$ over $B$'') is a (finite) sequence
$a_1, \ldots, a_\alpha \in \acl^n(A)$ together with a {\em core sequence of indices}
$0 = k_0 < k_1 < \ldots < k_m = \alpha$ such that
\begin{enumerate}
\item $\acl^n(A) \subseteq \acl^n(\{a_1, \ldots, a_\alpha\} \cup B)$ and,

\item for all $j = 0, \ldots, m-1$, 
$\{a_{k_j + 1}, \ldots, a_{k_{j+1}}\}$ is a basis of 
\[
\acl^n(A) \cap \{d \in M_n : \rk^n(d / \{a_1, \ldots, a_{k_j}\} \cup B) = 1 \}.
\]
\end{enumerate}
}\end{defin}

\begin{lem}\label{existence of n-ccs}
Let $A, B \subseteq M\meq$ where $A$ is finite.\\
(i) There is an $n$-ccs for $A/B$.\\
(ii) Suppose that $a_1, \ldots, a_\alpha$ and $a'_1, \ldots, a'_{\alpha'}$ are two $n$-ccs for $A/B$
with core sequences $0 = k_0 < k_1 < \ldots < k_s = \alpha$ and $0 = l_0 < l_1 < \ldots < l_t = \alpha'$,
respectively.
Then $s = t$ and, for all $i = 0, \ldots, s$, $k_i = l_i$ (so $\alpha = \alpha'$).
Moreover, for all $m = 1, \ldots, s$, 
$\acl\meq(\{a_1, \ldots, a_{k_m}\} \cup B) = \acl\meq(\{a'_1, \ldots, a'_{k_m}\} \cup B)$.
\end{lem}

\noindent
{\bf Proof.}
We prove part~(i) and point out during the argument why the uniqueness properties of part~(ii) follow.
Suppose that $\acl^n(A) \not\subseteq \acl^n(B)$ for otherwise the empty sequence has the required properties.
By convention let $k_0 = 0$.
By Lemma~\ref{finding a in acl of rank 1 over something else}
the following set is nonempty:
\[
X_1 = \acl^n(A) \cap \{d \in M_n : \rk^n(d / B) = 1 \}.
\]
By Lemma~\ref{n-rank 1 elements form a pregeometry}
$(\{d \in M_n : \rk^n(d / B) = 1 \}, \cl)$ with $\cl(Y) = \acl^n(YB) \cap \{d \in M_n : \rk^n(d / B) = 1 \}$ 
for all $Y \subseteq \{d \in M_n : \rk^n(d / B) = 1 \}$ is a pregeometry.
Let $k_1$ be the dimension of $X_1$ and let $\{a_1, \ldots, a_{k_1}\} \subseteq X_1$ be a basis of $X_1$.
Note that $k_1$ is determined only by $A$, $B$ and $n$.
Also observe that for any other basis  $\{a'_1, \ldots, a'_{k_1}\} \subseteq X_1$ we have 
$\acl\meq(\{a'_1, \ldots, a'_{k_1}\}  \cup B)  = \acl\meq(\{a_1, \ldots, a_{k_1}\}  \cup B)$.

Now suppose that, for some $m$, $k_1 < \ldots < k_m$ and $a_1, \ldots, a_{k_m} \in \acl^n(A)$
have been defined in such a way that, for all $j = 0, \ldots, m-1$,
$\{a_{k_j + 1}, \ldots, a_{k_{j+1}}\}$ is a basis of 
\[
X_j = \acl^n(A) \cap \{d \in M_n : \rk^n(d / \{a_1, \ldots, a_{k_j}\} \cup B) = 1 \}.
\]
Also suppose that the uniqueness properties of part~(ii) hold for $k_1, \ldots, k_m$ and 
$\acl\meq($ $\{a_1, \ldots, a_{k_i}\}$ $\cup B)$ for $i \leq m$.
If $\acl^n(A) \subseteq  \acl^n(\{a_1, \ldots, a_{k_m}\} \cup B)$ then we let $\alpha = k_m$ and
then $a_1, \ldots, a_\alpha$ is an $n$-ccs for $A/B$.

Suppose that $\acl^n(A) \not\subseteq \acl^n(\{a_1, \ldots, a_{k_m}\} \cup B)$.
By Lemma~\ref{finding a in acl of rank 1 over something else}
the set 
\[
X_{m+1} = \acl^n(A) \cap \{d \in M_n : \rk^n(d / \{a_1, \ldots, a_{k_m}\} \cup B) = 1 \} 
\]
is nonempty.
Let $d$ be the dimension of $X_{m+1}$, let $k_{m+1} = k_m + d$,  
and let 
\[
\{a_{k_m + 1}, \ldots, a_{k_{m+1}}\} \subseteq X_{m+1}
\] 
be a basis of $X_{m+1}$.
Note that $k_{m+1}$ depends only on $k_m$ and the dimension of $X_{m+1}$ where the latter depends only on 
$A$, $B$, $a_1, \ldots, a_{k_m}$ and $n$. 
The inductive assumption that for any choice of $a'_1, \ldots, a'_{k_m}$ with the same properties as 
$a_1, \ldots, a_{k_m}$ we have 
$\acl\meq(\{a'_1, \ldots, a'_{k_m}\} \cup B) = \acl\meq(\{a_1, \ldots, a_{k_m}\} \cup B)$
implies that 
$\acl\meq(\{a'_1, \ldots, a'_{k_{m+1}}\} \cup B) = \acl\meq(\{a_1, \ldots, a_{k_{m+1}}\} \cup B)$
for any basis $\{a'_{k_m + 1}, \ldots, a'_{k_{m+1}}\} \subseteq X_{m+1}$.
It also implies that, ultimately, $k_{m+1}$ depends only on $A$, $B$ and $n$.

Observe that for all $i = 1, \ldots, k_{m+1}$, 
$a_i \notin \acl^n(\{a_1, \ldots, a_{i-1}\} \cup B)$ so by 
Lemma~\ref{characterization of rank}
we have $\rk^n(A / B) \geq k_{m+1}$. 
Since $\rk^n(A / B)$ is finite (as $A$ is finite) it follows that the process of extending the
current sequence $a_1, \ldots, a_{k_{m+1}}$ in the described way will terminate after finitely many steps and
then we have an $n$-ccs for $A/B$.
\hfill $\square$

\begin{lem}\label{a long n-ccs is an n-cs}
Let $A, B \subseteq M\meq$ and suppose that $\rk^n(A / B) = \alpha < \omega$.
If $a_1, \ldots, a_\alpha$ is an $n$-ccs for $A/B$ then it is also an $n$-cs of $A/B$.
\end{lem}

\noindent
{\bf Proof.}
Suppose that $a_1, \ldots, a_\alpha$ is an $n$-ccs for $A/B$ with core sequence of indices
$0 = k_0 < k_1 < \ldots < k_m = \alpha$. 
From Definition~\ref{definition of n-ccs} of an $n$-ccs it follows that, for all $i = 1, \ldots, \alpha$,
$a_i \notin \acl^n(\{a_1, \ldots, a_{i-1}\} \cup B)$. 
Since $\alpha = \rk^n(A / B)$ it follows 
from Definition~\ref{definition of n-cs} 
of an $n$-cs that $a_1, \ldots, a_\alpha$ is an $n$-cs for $A/B$.
\hfill $\square$

\begin{lem}\label{for every n-cs there is an n-ccs}
Let $A, B \subseteq M\meq$ where $A$ is finite.\\
If $\rk^n(A / B) = \alpha$ then every $n$-ccs for $A/B$ has length $\alpha$.
\end{lem}

\noindent
{\bf Proof.}
Suppose that $\rk^n(A/B) = \alpha$.
By Lemma~\ref{existence of n-ccs},
it suffices to prove that there is at least one $n$-ccs for $A/B$ that has length $\alpha$.
By Lemma~\ref{characterization of rank} and
Definition~\ref{definition of n-cs} of $n$-cs,
there is an $n$-cs $a_1, \ldots, a_\alpha$ for $A/B$.
By the same definition and lemma we have 
\begin{align}\label{rank of a-i over a-i-1}
&a_i \notin \acl^n(\{a_1, \ldots, a_{i-1}\} \cup B) \ \text{ and } \\  
&\rk^n(a_i / \{a_1, \ldots, a_{i-1}\} \cup B) = 1
\ \ \text{ for all $i = 1, \ldots ,\alpha$}. \nonumber
\end{align}

\noindent
The rest of the proof will show how to transform $a_1, \ldots, a_\alpha$ into an $n$-ccs for $A/B$ with length $\alpha$.
This will be done step by step via a sequence of claims.

\begin{claim}\label{expanding the set by one elements}
Suppose that $0 \leq k_1 < k_2 < \alpha$ and that $\{a_{k_1 + 1}, \ldots, a_{k_2}\}$ is an independent
subset of 
\[
X = \acl^n(A) \cap \{d \in M_n : \rk^n(d / \{a_1, \ldots, a_{k_1}\} \cup B) = 1 \}.
\]
If $\{a_{k_1 + 1}, \ldots, a_{k_2}\}$ is {\em not} a basis of $X$ then there are $d \in X$
and $k_2 < l \leq \alpha$ such that $\{a_{k_1 + 1}, \ldots, a_{k_2}, d\}$ is an independent subset of $X$
and $a_1, \ldots, a_{l-1}, d, a_{l+1}, \ldots, a_\alpha$ is an $n$-cs for $A/B$.
\end{claim}

\noindent
{\em Proof.}
Suppose that $\{a_{k_1 + 1}, \ldots, a_{k_2}\}$ is {\em not} a basis of $X$.
There is $d \in X$ such that $d \notin \acl^n(\{a_1, \ldots, a_{k_2}\} \cup B)$.
If $d \notin \acl^n(\{a_1, \ldots, a_\alpha\} \cup B)$ then it follows 
(by~(\ref{rank of a-i over a-i-1}) and Lemma~\ref{characterization of rank})
that $\rk^n(A / B) \geq \alpha + 1$, which contradicts the assumption.
Hence there is a minimal $1 \leq l \leq \alpha$ such that
$d \in \acl^n(\{a_1, \ldots, a_l\} \cup B)$.
Since $d \notin \acl^n(\{a_1, \ldots, a_{k_2}\} \cup B)$ we must have $l > k_2$.
Then $d \notin \acl^n(\{a_1, \ldots, a_{l-1}\} \cup B)$.
As also $d \in X$ it follows (using Lemma~\ref{monotonicity of n-rank}) that
$\rk^n(d / \{a_1, \ldots, a_{l-1}\} \cup B) = 1$.
By~(\ref{rank of a-i over a-i-1}), we have
$\rk^n(a_l / \{a_1, \ldots, a_{l-1}\} \cup B) = 1$.
Now Assumption~\ref{assumption about exchange property}
implies that $a_l \in \acl^n( \{a_1, \ldots, a_{l-1}, d\} \cup B)$.

Suppose, for a contradiction, that there is $i > l$ such that 
\[
a_i \in \acl^n(\{a_1, \ldots, a_{l-1}, d, a_{l+1}, \ldots, a_{i-1}\} \cup B).
\]
Since $d \in \acl^n(\{a_1, \ldots, a_l\} \cup B)$ it follows that
$a_i \in \acl^n(\{a_1, \ldots, a_{i-1}\} \cup B)$ which contradicts~(\ref{rank of a-i over a-i-1}).
Hence, no such $i > l$ exists and it follows that 
$a_1, \ldots, a_{l-1}, d, a_{l+1}, \ldots, a_\alpha$ is an $n$-cs for $A/B$.
\hfill $\square$

\begin{claim}\label{reordering one element}
Suppose that, for some $0 \leq k < \alpha$ and $1 \leq s \leq \alpha - k$, we have 
\[
\rk^n(a_{k+i} / \{a_1, \ldots, a_k\} \cup B) = 1 \ \text{ for all $i = 1, \ldots, s$.}
\]
If $k+s < l \leq \alpha$ and $\rk^n(a_l / \{a_1, \ldots, a_k\} \cup B) = 1$
then 
\[
a_1, \ldots, a_{k+s}, a_l, a_{k+s+1}, \ldots, a_{l-1}, a_{l+1}, \ldots, a_\alpha
\ \ \text{ is an $n$-cs for $A/B$.}
\]
\end{claim}

\noindent
{\em Proof.}
Let $k$ and $s$ be as assumed.
Suppose that $k+s < l \leq \alpha$ and $\rk^n(a_l / \{a_1, \ldots, a_k\} \cup B) = 1$.
Suppose, for a contradiction, that there is $k+s < j < l$ such that
$a_j \in \acl^n(\{a_1, \ldots, a_{k+s}, a_l, a_{k+s+1}, \ldots, a_{j-1}\} \cup B)$.
By~(\ref{rank of a-i over a-i-1}) and Lemma~\ref{monotonicity of n-rank}, 
\begin{align*}
1 &= \rk^n(a_l / \{a_1, \ldots, a_k\} \cup B) \\
&\geq 
\rk^n(a_l / \{a_1, \ldots, a_{j-1}\} \cup B) \\
&\geq 
\rk^n(a_l / \{a_1, \ldots, a_{l-1}\} \cup B) = 1.
\end{align*}
Hence $\rk^n(a_l / \{a_1, \ldots, a_{j-1}\} \cup B) = 1$.
By~(\ref{rank of a-i over a-i-1}) we also have
$\rk^n(a_j / \{a_1, \ldots, a_{j-1}\} \cup B) = 1$.
Now Assumption~\ref{assumption about exchange property}
implies that $a_l \in \acl^n(\{a_1, \ldots, a_j\} \cup B)$ and since $j < l$ this 
contradicts~(\ref{rank of a-i over a-i-1}).

For a contradiction, suppose that there is $l < j \leq \alpha$ such that
\[
a_j \in \acl^n(\{a_1, \ldots, a_{k+s}, a_l, a_{k+s+1}, \ldots, a_{l-1}, a_{l+1}, \ldots, a_{j-1}\} \cup B).
\]
Since $\{a_1, \ldots, a_{k+s}, a_l, a_{k+s+1}, \ldots, a_{l-1}, a_{l+1}, \ldots, a_{j-1}\} = \{a_1, \ldots a_{j-1}\}$
we have a contradiction to~(\ref{rank of a-i over a-i-1}).
It follows that $a_1, \ldots, a_{k+s}, a_l, a_{k+s+1}, \ldots, a_{l-1}, a_{l+1}, \ldots, a_\alpha$ 
is an $n$-cs for $A/B$.
\hfill $\square$

\begin{claim}\label{expanding the set by one element and reordering}
Suppose that $0 \leq k_1 < k_2 < \alpha$ and that $\{a_{k_1 + 1}, \ldots, a_{k_2}\}$ is an independent
subset of 
\[
X = \acl^n(A) \cap \{d \in M_n : \rk^n(d / \{a_1, \ldots, a_{k_1}\} \cup B) = 1\}.
\]
If $\{a_{k_1 + 1}, \ldots, a_{k_2}\}$ is {\em not} a basis of $X$ then there is $d \in X$ such that
$\{a_{k_1 + 1}, \ldots, a_{k_2}, d\}$ is an independent subset of $X$ and, for some $k_2 < i \leq \alpha$,
\[
a_1, \ldots, a_{k_2}, d, a_{k_2 + 1}, \ldots, a_{i-1}, a_{i+1}, \ldots, a_\alpha
\ \ \text{ is an $n$-cs for $A/B$.}
\]
\end{claim}

\noindent
{\em Proof.}
Let us make the assumptions of the claim.
Moreover, suppose that $\{a_{k_1 + 1}, \ldots, a_{k_2}\}$ is not a basis of $X$.
Now we first apply Claim~\ref{expanding the set by one elements}
to find $d \in X$
and $k_2 < i \leq \alpha$ such that $\{a_{k_1 + 1}, \ldots, a_{k_2}, d\}$ is an independent subset of $X$
and $a_1, \ldots, a_{i-1}, d, a_{i+1}, \ldots, a_\alpha$ is an $n$-cs for $A/B$.
Then Claim~\ref{reordering one element} implies that
also $a_1, \ldots, a_{k_2}, d, a_{k_2 + 1}, \ldots, a_{i-1}, a_{i+1}, \ldots, a_\alpha$
is an $n$-cs for $A/B$.
\hfill $\square$

\begin{claim}\label{expanding the set by several elements}
Suppose that $0 \leq k_1 < k_2 < \alpha$ and that $\{a_{k_1 + 1}, \ldots, a_{k_2}\}$ is an independent
subset of 
\[
X = \acl^n(A) \cap \{d \in M_n : \rk^n(d / \{a_1, \ldots, a_{k_1}\} \cup B) = 1\}.
\]
If $\{a_{k_1 + 1}, \ldots, a_{k_2}\}$ is {\em not} a basis of $X$ then there are $1 \leq s \leq \alpha - k_2$,
$d_1, \ldots, d_s \in X$, and $a'_{\alpha - k_2 - s}, \ldots, a'_\alpha \in M_n$ such that 
$\{a_{k_1 + 1}, \ldots, a_{k_2}, d_1, \ldots, d_s\}$ is a basis of $X$ and
$a_1, \ldots, a_{k_2}, d_1, \ldots, d_s, a'_{\alpha - k_2 - s}, \ldots, a'_\alpha$
is an $n$-cs for $A/B$.
\end{claim}

\noindent
{\em Proof.}
The claim follows from repeated uses of 
Claim~\ref{expanding the set by one element and reordering}
and the fact that all bases of $X$ have the same finite cardinality.
\hfill $\square$

\begin{claim}\label{inductive step in finding an n-ccs}
Suppose that there are $m < \alpha$ and  $0 = k_0 < k_1 < \ldots < k_m < \alpha$ such that, for all $j = 0, \ldots, m-1$,
$\{a_{k_j + 1}, \ldots, a_{k_{j+1}}\}$ is a basis of 
\[
X_j = \acl^n(A) \cap \{ d \in M_n : \rk^n(d / \{a_1, \ldots, a_{k_j}\} \cup B) = 1\}.
\]
Then there are $k_m < k_{m+1} \leq \alpha$, a basis $\{a'_{k_m + 1}, \ldots, a'_{k_{m+1}}\}$ of 
\[
X_{m+1} =  \acl^n(A) \cap \{ d \in M_n : \rk^n(d / \{a_1, \ldots, a_{k_m}\} \cup B) = 1\},
\]
and $a'_{k_{m+1} + 1}, \ldots, a'_\alpha \in M_n$ such that
$a_1, \ldots, a_{k_m}, a'_{k_m + 1}, \ldots, a'_\alpha$ is an $n$-cs for $A/B$.
\end{claim}

\noindent
{\em Proof.}
From~(\ref{rank of a-i over a-i-1}) we get
$\rk^n(a_{k_m + 1} / \{a_1, \ldots, a_{k_m}\} \cup B) = 1$.
If $\{a_{k_m + 1}\}$ is a basis of $X_{m+1}$ then let $a'_i = a_i$ for all $i = k_m + 1, \ldots, \alpha$
and we are done.
Otherwise we get the conclusion of the claim by 
Claim~\ref{expanding the set by several elements}.
\hfill $\square$

\medskip

\noindent
Observe that the assumptions of Claim~\ref{inductive step in finding an n-ccs}
are vacuously satisfied if $m = 0$.
Therefore Lemma~\ref{for every n-cs there is an n-ccs}
follows by induction on $m$ where Claim~\ref{inductive step in finding an n-ccs}
serves as the inductive step. 
This completes the proof of Lemma~\ref{for every n-cs there is an n-ccs}.
\hfill $\square$

\begin{lem}\label{n-ccs and rk-n}
Let $A, B, C \subseteq M\meq$ where $B$ is finite.
Suppose that $\rk^n(B / AC) = \rk^n(B / C) = \beta$ (so $\beta < \omega$).\\
(i) If $b_1, \ldots, b_\beta$ is an $n$-ccs for $B/AC$ with core sequence of indices
$0 = k_0 < k_1 < \ldots < k_t = \beta$, then $b_1, \ldots, b_\beta$ is also an $n$-ccs for $B/C$
with the same core sequence of indices.\\
(ii) If $b_1, \ldots, b_\beta$ is an $n$-ccs for $B/C$ with core sequence of indices
$0 = k_0 < k_1 < \ldots < k_t = \beta$, then $b_1, \ldots, b_\beta$ is also an $n$-ccs for $B/AC$
with the same core sequence of indices.
\end{lem}

\noindent
{\bf Proof.}
Suppose that $B$ is finite and $\rk^n(B / AC) = \rk^n(B / C) = \beta$.

(i) Let $b_1, \ldots, b_\beta$ be an $n$-ccs for $B/AC$ with core sequence of indices
$0 = k_0 < k_1 < \ldots < k_t = \beta$.
Then, for all $m = 0, \ldots, t-1$, 
\begin{align}\label{defining properties of the n-ccs b}
&\{b_{k_m + 1}, \ldots, b_{k_{m+1}}\} \ \text{ is a basis of } \\
&X_m^{AC} = \acl^n(B) \cap \{d \in M_n : \rk^n(d / \{b_1, \ldots, b_{k_m}\} \cup AC) = 1 \} \ \text{ and,}   \nonumber \\
&\text{for all } i = k_m + 1, \ldots, k_{m+1}, \rk^n(b_i / \{b_1, \ldots, b_{k_m}\} \cup AC) = 1. \nonumber
\end{align}
This implies that
\begin{align}\label{b-i is not in the closure of b-i-1 union AC}
&\text{for all } \ i = 1, \ldots, \beta, b_i \notin \acl^n(\{b_1, \ldots, b_{i-1}\} \cup AC) \ \text{ and hence} \\
&b_i \notin \acl^n(\{b_1, \ldots, b_{i-1}\} \cup C). \nonumber
\end{align}
It also follows (using Lemma~\ref{monotonicity of n-rank})
that,  for all $m = 0, \ldots, t-1$ and all $ i = k_m + 1, \ldots, k_{m+1}$,
\[
\rk^n(b_i / \{b_1, \ldots, b_{k_m}\} \cup C) \geq 1.
\]
Towards a contradiction, suppose that there are $m \in \{0, \ldots, t-1\}$ and $i \in \{k_m +1, \ldots, k_{m+1}\}$
such that $\rk^n(b_i / \{b_1, \ldots, b_{k_m}\} \cup C) \geq 2$.
Then there is $d \in \acl^n(b_i)$ such that $d \notin \acl^n(\{b_1, \ldots, b_{k_m}\} \cup C)$ and
$\rk^n(b_i / \{d, b_1, \ldots, b_{k_m}\} \cup C) \geq 1$. 
Since $d \in \acl^n(b_i)$ it follows from~(\ref{b-i is not in the closure of b-i-1 union AC}) that
for all $k_m < j \leq \alpha$, 
\[
b_j \notin \acl^n(\{b_1, \ldots, b_{k_m}, d, b_{k_m + 1}, \ldots, b_{j-1}\} \cup C).
\]
This implies (via Lemma~\ref{characterization of rank}) that $\rk^n(B / C) \geq \beta + 1$,
contradicting the assumption.
Thus we conclude that
\begin{align}\label{all b-i have rank 1 over C}
&\text{for all $m =  0, \ldots, t-1$ and all $ i = k_m + 1, \ldots, k_{m+1}$,} \\
&\rk^n(b_i / \{b_1, \ldots, b_{k_m}\} \cup C) = 1. \nonumber
\end{align}
By~(\ref{b-i is not in the closure of b-i-1 union AC}) and~(\ref{all b-i have rank 1 over C}),
for all $m =  0, \ldots, t-1$,
$\{b_{k_m + 1}, \ldots, b_{k_{m+1}}\}$ is an independent subset of 
\[
X^C_m = \acl^n(B) \cap \{d \in M_n : \rk^n(d / \{b_1, \ldots, b_{k_m}\} \cup C) = 1\}.
\]

Towards a contradiction, suppose that for some $m \in \{0, \ldots, t-1 \}$
there is $d \in X^C_m$ such that $d \notin \acl^n(\{b_1, \ldots, b_{k_m}\} \cup C)$.
Since $\{b_1, \ldots, b_{k_m}\}$ is a basis of $X_m^{AC}$ it follows that $d \in \acl^n(\{b_1, \ldots, b_{k_m}\} \cup AC)$.
If, for some $i > k_m$, we would have that
$b_i \in \acl^n(\{b_1, \ldots, b_{k_m}, d, b_{k_m + 1}, \ldots, b_{i-1}\} \cup C)$, 
then we would get 
\[
b_i \in \acl^n(\{b_1, \ldots, b_{k_m}, b_{k_m + 1}, \ldots, b_{i-1}\} \cup AC)
\]
which would contradict~(\ref{b-i is not in the closure of b-i-1 union AC}).
Hence we conclude that for all $i > k_m$,
\[
b_i \notin \acl^n(\{b_1, \ldots, b_{k_m}, d, b_{k_m + 1}, \ldots, b_{i-1}\} \cup C).
\]
But this implies that $\rk^n(B / C) \geq \beta + 1$ which contradicts the assumption.
Hence there cannot be any  $m \in \{0, \ldots, t-1 \}$
and $d \in X^C_m$ such that $d \notin \acl^n(\{b_1, \ldots, b_{k_m}\} \cup C)$,
and therefore $\{b_1, \ldots, b_{k_m}\}$ is a basis of $X^C_m$.
It follows that $b_1, \ldots, b_\beta$ is an $n$-ccs for $B/C$ with  core sequence of indices
$0 = k_0 < k_1 < \ldots < k_t = \beta$.

(ii) Now suppose that $b_1, \ldots, b_\beta$ is an $n$-ccs for $B/C$ with core sequence of indices
$0 = k_0 < k_1 < \ldots < k_\delta = \beta$.
It follows that for all $m = 0, \ldots, \delta - 1$ and all $i = k_m +1, \ldots, k_{m+1}$,
\begin{align}\label{properties of the b over B/C}
&\rk^n(b_i / \{b_1, \ldots, b_{i-1}\} \cup C) = \rk^n(b_i / \{b_1, \ldots, b_{k_m}\} \cup C) = 1, \ \text{ and } \\
&b_i \notin \acl^n(\{b_1, \ldots, b_{i-1}\} \cup C). \nonumber
\end{align}
By Lemmas~\ref{existence of n-ccs}
and~\ref{for every n-cs there is an n-ccs},
there is an $n$-ccs $b'_1, \ldots, b'_\beta$ for $B/AC$ with core sequence of indices
$0 = l_0 < l_1 < \ldots < l_\gamma = \beta$.
By part~(i), $b'_1, \ldots, b'_\beta$ is also an $n$-ccs for $B/C$ with the same core sequence of indices
$0 = l_0 < l_1 < \ldots < l_\gamma = \beta$.
By Lemma~\ref{existence of n-ccs},
we then have $\gamma = \delta$ and $l_m = k_m$ for all $m = 0, \ldots, \delta$.
The same lemma also tells that, for all $m = 1, \ldots, \delta$, 
\begin{equation}\label{acl of b up to k-m and C}
\acl\meq(\{b_1, \ldots, b_{k_m}\} \cup C) = \acl\meq(\{b'_1, \ldots, b'_{k_m}\} \cup C).
\end{equation}
Now we prove that 
\begin{equation}\label{acl of b up to k-m and AC}
\acl\meq(\{b_1, \ldots, b_{k_m}\} \cup AC) = \acl\meq(\{b'_1, \ldots, b'_{k_m}\} \cup AC).
\end{equation}
Suppose that $d \in \acl\meq(\{b_1, \ldots, b_{k_m}\} \cup AC)$.
Then there are a formula $\varphi(u, \bar{x}, \bar{y}, \bar{z})$,
a tuple $\bar{b}$ of elements from $\{b_1, \ldots, b_{k_m}\}$,
a tuple $\bar{a}$ of elements from $A$, and a tuple $\bar{c}$ of elements from $C$,
such that $\mcM\meq \models \varphi(d, \bar{b}, \bar{a}, \bar{c})$ and
only finitely many elements satisfy $\varphi(u, \bar{b}, \bar{a}, \bar{c})$.
Let $m$ be such that $d$ and all elements that belong to any of $\bar{b}$, $\bar{a}$, or $\bar{c}$,
are members of $M_m$. As $\mcM_m$ is $\omega$-categorical there is a formula
$\theta(u, \bar{x}, \bar{y}, \bar{z})$ which isolates $\tp_{\mcM_m}(d, \bar{b}, \bar{a}, \bar{c})$,
where `$\tp_{\mcM_m}$' denotes the type computed in $\mcM_m$.
Then only finitely many elements satisfy $\theta(u, \bar{b}, \bar{a}, \bar{c})$,
all of them belong to $M_m$, and $d$ is one among them.
From~(\ref{acl of b up to k-m and C})
we have $\rng(\bar{b}) \subseteq \acl^m(\{b'_1, \ldots, b'_{k_m}\} \cup C)$, so 
there are a formula $\psi(\bar{x}, \bar{y}, \bar{z})$, a tuple $\bar{b}'$
of elements from $\{b'_1, \ldots, b'_{k_m}\}$, and a tuple $\bar{c}'$ of elements from $C$ such
that $\mcM\meq \models \psi(\bar{b}, \bar{b}', \bar{c}')$ and $\psi(\bar{x}, \bar{b}', \bar{c}')$
is satisfied by only finitely many tuples.
It is now clear that $d$ satisfies the formula
\[
\chi(u, \bar{a}, \bar{b}', \bar{c}, \bar{c}') :=
\exists \bar{x} \big[ \theta(u, \bar{x}, \bar{a}, \bar{c}) \wedge \psi(\bar{x}, \bar{b}', \bar{c}')\big].
\]
The formula $\psi(\bar{x}, \bar{b}', \bar{c}')$ is satisfied by only finitely many tuples,
and if $\bar{b}^*$ is one of them and $\mcM\meq \models  \theta(d', \bar{b}^*, \bar{a}, \bar{c})$
for some $d'$, then $d'$ and all members of $\bar{b}^*$ belong to $M_m$, and
$\tp_{\mcM_m}(d', \bar{b}^*, \bar{a}, \bar{c}) = \tp_{\mcM_m}(d, \bar{b}, \bar{a}, \bar{c})$,
so only finitely many elements satisfy $\theta(u, \bar{b}^*, \bar{a}, \bar{c})$.
It now follows that $\chi(u, \bar{a}, \bar{b}', \bar{c}, \bar{c}')$ is satisfied by only finitely many elements.
Hence 
\[
d \in \acl\meq(\bar{a}, \bar{b}', \bar{c}, \bar{c}') \subseteq 
\acl\meq(\{b'_1, \ldots, b'_{k_m}\} \cup AC).
\]
This shows the inclusion from left to right in~(\ref{acl of b up to k-m and AC})
and by a similar argument (letting the $b_i$ and $b'_i$ switch roles) the converse inclusion follows.
This proves~(\ref{acl of b up to k-m and AC}).

For all $m = 0, \ldots, \delta - 1$, let 
\[
X_m^{AC} = \acl^n(B) \cap \{d \in M_n : \rk^n(d / \{b_1, \ldots, b_{k_m}\} \cup AC) = 1\}
\]
and note that by~(\ref{acl of b up to k-m and AC}) we also have
\begin{equation}\label{X-m-AC can be defined in terms of b'}
X_m^{AC} = \acl^n(B) \cap \{d \in M_n : \rk^n(d / \{b'_1, \ldots, b'_{k_m}\} \cup AC) = 1\}.
\end{equation}
Since $b'_1, \ldots, b'_\beta$ is an $n$-ccs for $B/AC$ with core sequence of indices
$0 = k_0 < k_1 < \ldots < k_\delta = \beta$
it follows that, for all $m = 0, \ldots, \delta - 1$, 
$\{b'_{k_m + 1}, \ldots, b'_{k_{m+1}}\}$ is a basis of $X_m^{AC}$, so $X_m^{AC}$ has dimension
$k_{m+1} - k_m$.
From~(\ref{acl of b up to k-m and AC})
it follows that 
\[
X_m^{AC} \subseteq \acl^n(\{b_1, \ldots, b_{k_{m+1}}\} \cup AC).
\]
Fix any $m \in \{0, \ldots, \delta - 1\}$.
From~(\ref{properties of the b over B/C}), we get
\[
\rk^n(b_i / \{b_1, \ldots, b_{k_m}\} \cup AC) \leq 1 \ \text{ for all $i = k_m + 1, \ldots, k_{m+1}$}.
\]
Suppose, for a contradiction, that there is  $i \in \{ k_m + 1, \ldots, k_{m+1}\}$ such that
\begin{itemize}
\item $\rk^n(b_i / \{b_1, \ldots, b_{k_m}\} \cup AC) = 0$
(equivalently, $b_i \in \acl^n(\{b_1, \ldots, b_{k_m}\} \cup AC)$), 
\item {\bf or} $b_i \in \acl^n(B_iAC)$ where $B_i = \{b_{k_m + 1}, \ldots, b_{k_{m+1}}\} \setminus \{b_i\}$.
\end{itemize}
Then some proper subset of  $\{b_{k_m + 1}, \ldots, b_{k_{m+1}}\}$ is a basis of $X_m^{AC}$
which contradicts our previous conclusion that $X_m^{AC}$ has dimension $k_{m+1} - k_m$.
Hence we conclude that $\rk^n(b_i / \{b_1, \ldots, b_{k_m}\} \cup AC) = 1$ for all $i = k_m + 1, \ldots, k_{m+1}$
and that $\{b_{k_m + 1}, \ldots, b_{k_{m+1}}\}$ is independent 
(in the pregeometry $\{d \in M_n : \rk^n(d / \{b_1, \ldots, b_{k_m}\} \cup AC) = 1\}$)
and hence a basis of $X_m^{AC}$.
Since the argument holds for all $m = 0, \ldots, \delta - 1$ it follows that
$b_1, \ldots, b_\beta$ is an $n$-ccs for $B/AC$ with core sequence of indices 
$0 = k_0 < k_1 < \ldots < k_\delta = \beta$.
\hfill $\square$

\begin{prop}\label{symmetry of n-rank}
Let $A, C \subseteq M\meq$ and $B \subseteq M_n$ where $A$ and $B$ are finite.
If $\rk^n(A / BC) < \rk^n(A / C)$ then $\rk^n(B / AC) < \rk^n(B / C)$.
\end{prop}

\noindent
{\bf Proof.}
Let $A, C \subseteq M\meq$ and $B \subseteq M_n$ where $A$ and $B$ are finite.
We prove the result by induction on $\rk^n(B / C)$.
If $\rk^n(B / C) = 0$ then, since $B  \subseteq M_n$ and by the definition of $\rk^n$, 
$B \subseteq \acl^n(B) \subseteq \acl^n(C)$ 
so $\acl^n(BC) = \acl^n(C)$. Hence 
(by the definition of $\rk^n$)
$\rk^n(A / BC) = \rk^n(A / C)$ and the statement is vacuously satisfied.

The induction step remains. The induction hypothesis will be:
\begin{enumerate}
\item[(IH)] For every finite $B' \subseteq M_n$, if $\rk^n(B' / C) < \rk^n(B / C)$ and $\rk^n(A / B'C) < \rk^n(A / C)$,
then $\rk^n(B' / AC) < \rk^n(B' / C)$.
\end{enumerate}
Suppose that $\rk^n(A / BC) < \rk^n(A / C)$.
Also let $\beta = \rk^n(B / C)$. (We will show that $\rk^n(B / AC) < \beta$.)
Then there is an $n$-ccs $b_1, \ldots, b_\beta$ for $B/C$ with core sequence of indices
$0 = k_0 < k_1 < \ldots < k_\delta = \beta$.
In particular, $b_1, \ldots, b_\beta \in \acl^n(B)$.
We divide the argument into two main cases.

\medskip
\noindent
{\bf Case 1.} Suppose that there is $s < \beta$ such that 
$\rk^n(A / \{b_1, \ldots, b_s\} \cup C) < \rk^n(A / C)$.

\medskip

\noindent
We may assume that $s$ is minimal such that the above holds.
As $b_1, \ldots, b_\beta$ is an $n$-ccs for $B/C$, hence
(by Lemma~\ref{a long n-ccs is an n-cs}) 
an $n$-cs for $B/C$,
it follows from Lemma~\ref{characterization of rank}
that
\[
\rk^n(b_1, \ldots, b_s / C) = s < \beta = \rk^n(B / C).
\]
Now the induction hypothesis (IH) implies that
\[
\rk^n(b_1, \ldots, b_s / AC) < \rk^n(b_1, \ldots, b_s / C) = s.
\]
By Lemma~\ref{characterization of rank}
there is $t \leq s$ such that 
\begin{equation}\label{b-t belongs to acl of b-t-1 AC}
b_t \in \acl^n(\{b_1, \ldots, b_{t-1}\} \cup AC).
\end{equation}
Let $t$ be minimal such that the above holds.
Towards a contradiction, suppose that $\rk^n(B / AC) = \rk^n(B / C) = \beta$.
Then there is an $n$-ccs $b'_1, \ldots, b'_\beta$ for $B/AC$ with core sequence of indices
$0 = l_0 < l_1 < \ldots < l_\gamma = \beta$.
By Lemma~\ref{n-ccs and rk-n},
$b_1, \ldots, b_\beta$ is also an $n$-ccs for $B/C$ with core sequence of indices
$0 = k_0 < k_1 < \ldots < k_\gamma = \beta$.
By Lemma~\ref{existence of n-ccs}
we have $\gamma = \delta$ and $l_i = k_i$ for all $i = 0, \ldots, \delta$.
From Lemma~\ref{existence of n-ccs} it follows that
if, for any $m \in \{0, \ldots, \delta - 1\}$, we define 
\[
X_m^{AC} = \acl^n(B) \cap \{d \in M_n : \rk^n(d / \{b_1, \ldots, b_{k_m}\} \cup AC) = 1 \}
\]
then
\[
X_m^{AC} = \acl^n(B) \cap \{d \in M_n : \rk^n(d / \{b'_1, \ldots, b'_{k_m}\} \cup AC) = 1 \}.
\]
Also, $b_{k_m + 1}, \ldots, b_{k_{m+1}}, b'_{k_m + 1}, \ldots, b'_{k_{m+1}} \in X_m^{AC}$.
Let $m$ be such that $k_m < t \leq k_{m+1}$.
Now~(\ref{b-t belongs to acl of b-t-1 AC}) implies that
some {\em proper} subset of $\{b_{k_m + 1}, \ldots, b_{k_{m+1}}\}$ is a basis of $X_m^{AC}$.
But since $b'_1, \ldots, b'_\beta$ is an $n$-ccs for $B/AC$ with core sequence of indices
$0 = k_0 < k_1 < \ldots < k_\delta = \beta$ it follows that
$\{b'_{k_m + 1}, \ldots, b'_{k_{m+1}}\}$ is a basis of $X_m^{AC}$.
So there are two bases of $X_m^{AC}$ with different cardinalities, which is impossible.
Thus we conclude that $\rk^n(B / AC) < \beta$.

\medskip
\noindent
{\bf Case 2.} Suppose that Case 1 does not hold, that is, suppose that for all $s < \beta$,
$\rk^n(A / \{b_1, \ldots, b_s\} \cup C) = \rk^n(A / C)$.

\medskip

\noindent
Let $\alpha = \rk^n(A / C)$, so by assumption,
\[
rk^n(A / \{b_1, \ldots, b_{\beta - 1}\} \cup C) = \rk^n(A / C) = \alpha.
\]
Since $b_1, \ldots, b_\beta$ is an $n$-ccs for $B/C$, hence an $n$-cs for $B/C$, it follows 
(by Lemma~\ref{characterization of rank} and the assumption that $B\subseteq M_n$) 
that 
$B \subseteq \acl^n(B) \subseteq \acl^n(\{b_1, \ldots, b_\beta\} \cup C)$ and hence
$\rk^n(A / BC) = \rk^n(A / \{b_1, \ldots, b_\beta\} \cup C)$.
From the assumption that $\rk^n(A / BC) < \rk^n(A /C)$ we now get
\begin{equation}\label{rank over b-beta C is less than over b-beta-1 C}
\rk^n(A / \{b_1, \ldots, b_\beta\} \cup C) < \rk^n(A / \{b_1, \ldots, b_{\beta-1}\} \cup C) = \alpha.
\end{equation}
Let $a_1, \ldots, a_\alpha$ be an $n$-ccs for $A/C$ with core sequence of indices
$0 = l_0 < l_1 < \ldots < l_\gamma = \alpha$.
By Lemma~\ref{n-ccs and rk-n},
$a_1, \ldots, a_\alpha$ is also an $n$-ccs for $A/\{b_1, \ldots, b_{\beta-1}\}\cup C$
with the same core sequence of indices 
$0 = l_0 < l_1 < \ldots < l_\gamma = \alpha$.
This implies that,
\[
\text{ for all $i = 1, \ldots, \alpha$, } \ \rk^n(a_i / \{a_1, \ldots, a_{i-1}, b_1, \ldots, b_{\beta-1}\} \cup C) = 1.
\]
Since $b_1, \ldots, b_\beta$ is an $n$-ccs for $B/C$ also have that
\begin{equation}\label{b-i over earlier b and C have rank 1, another time perhaps}
\text{ for all $i = 1, \ldots, \beta$, } \ \rk^n(b_i / \{b_1, \ldots, b_{i-1}\} \cup C) = 1.
\end{equation}
By~(\ref{rank over b-beta C is less than over b-beta-1 C})
(and Lemma~\ref{characterization of rank}) 
there is $s \leq \alpha$ such that 
\begin{equation}\label{a-s is in the acl of previous a, b and C}
a_s \in \acl^n(\{a_1, \ldots, a_{s-1}, b_1, \ldots, b_\beta\} \cup C).
\end{equation}
Let $s$ be minimal such that the above holds.

If $b_\beta \in \acl^n(\{b_1, \ldots, b_{\beta-1}\} \cup AC)$ then $b_1, \ldots, b_\beta$ 
is not an $n$-ccs for $B/AC$ so by 
Lemma~\ref{n-ccs and rk-n}
we get $\rk^n(B / AC) < \rk^n(B / C)$.

Now suppose that $b_\beta \notin \acl^n(\{b_1, \ldots, b_{\beta-1}\} \cup AC)$.
As $a_1, \ldots, a_{s-1} \in \acl^n(A)$ we get 
$b_\beta \notin \acl^n(\{b_1, \ldots, b_{\beta-1}, a_1, \ldots, a_{s-1}\} \cup C)$ and hence
$\rk^n(b_\beta / \{b_1, \ldots, b_{\beta-1}, a_1, \ldots, a_{s-1}\} \cup C) \geq 1$.
This and~(\ref{b-i over earlier b and C have rank 1, another time perhaps})
(together with Lemma~\ref{monotonicity of n-rank})
gives
\begin{equation}\label{rank of b-beta over b and a and C}
\rk^n(b_\beta / \{b_1, \ldots, b_{\beta-1}, a_1, \ldots, a_{s-1}\} \cup C) = 1.
\end{equation}
Since $a_1, \ldots, a_\alpha$ is an $n$-ccs for $A/\{b_1, \ldots, b_{\beta-1}\} \cup C$ 
we get
\begin{equation}\label{rank of a-s over b and a and C}
\rk^n(a_s / \{b_1, \ldots, b_{\beta-1}, a_1, \ldots, a_{s-1}\} \cup C) = 1.
\end{equation}
Now (\ref{a-s is in the acl of previous a, b and C}), (\ref{rank of b-beta over b and a and C}),
(\ref{rank of a-s over b and a and C}) and 
Assumption~\ref{assumption about exchange property}
imply that
\[
b_\beta \in \acl^n(\{b_1, \ldots, b_{\beta-1}, a_1, \ldots, a_s\} \cup C)
\]
and therefore $b_\beta \in \acl^n(\{b_1, \ldots, b_{\beta-1}\} \cup AC)$.
Then $b_1, \ldots, b_\beta$ is not an $n$-ccs for $B/AC$, so by 
Lemma~\ref{n-ccs and rk-n} it follows that
$\rk^n(B / AC) < \rk^n(B / C)$.
This completes the proof of Proposition~\ref{symmetry of n-rank}.
\hfill $\square$

\begin{prop}\label{symmetry of n-independence} {\bf (Symmetry of $\indn$ restricted to $M_n$)}
Let $A, C \subseteq M\meq$ and $B \subseteq M_n$.
If $A \underset{C}{\nindn} B$ then $B \underset{C}{\nindn} A$.
\end{prop}

\noindent
{\bf Proof.}
Suppose that $A \underset{C}{\nindn} B$ where $B \subseteq M_n$.
Then there is finite $A' \subseteq A$ such that $\rk^n(A' / BC) < \rk^n(A' / C)$.
By Lemma~\ref{locality property of rank},
there is finite $B' \subseteq B$ such that $\rk^n(A' / B'C) < \rk^n(A' / C)$.
By Lemma~\ref{monotonicity of n-rank} and Proposition~\ref{symmetry of n-rank},
\[
\rk^n(B' / AC) \leq \rk^n(B' / A'C) < \rk^n(B' / C)
\] 
and hence 
$B \underset{C}{\nindn} A$.
\hfill $\square$

\medskip
\noindent
In Section~\ref{n-independence} we saw that, 
even without Assumption~\ref{assumption about exchange property},
$\indn$ has all the properties of an independence relation restricted to $M_n$
{\em except}, possibly, for the symmetry property.
Now we have:

\begin{theor}\label{n-independence is an independence relation}
Let $T$ be $\omega$-categorical, let $n < \omega$, and suppose that 
Assumption~\ref{assumption about exchange property} holds for this $n$.
Then $\indn$ is an independence relation restricted to $M_n$.
\end{theor}

\noindent
{\bf Proof.}
By Proposition~\ref{symmetry of n-independence}, $\indn$ has the symmetry property with respect to
subsets of $M_n$.
By the results in Section~\ref{n-independence}
$\indn$ has all the other properties of an independence relation, with respect to any subsets of $M\meq$.
\hfill $\square$

\medskip
\noindent
Under Assumption~\ref{assumption about exchange property} we can strengthen 
one part of Lemma~\ref{monotonicity of n-rank} as follows:

\begin{lem}\label{modularity of rank}
Suppose that $A, B, C \subseteq M\meq$ where $A$ is finite and $B \subseteq A$.
Then $\rk^n(A / C) = \rk^n(B / C) + \rk^n(A / BC)$.
\end{lem}

\noindent
{\bf Proof.}
Let $\beta = \rk^n(B / C)$ and $\alpha = \rk^n(A / BC)$.
By Lemmas~\ref{existence of n-ccs} and~\ref{for every n-cs there is an n-ccs}
there is an $n$-ccs $b_1, \ldots, b_\beta$ for $B/C$ with
core sequence of indices $0 = k_0 < k_1 < \ldots < k_s = \beta$.
By the same lemmas there is also an $n$-ccs $a_1, \ldots, a_\alpha$ for $A/BC$ with core sequence
of indices $0 = l_0 < l_1 < \ldots < l_t = \alpha$.
Since we must have $\acl^n(B) \subseteq \acl^n(\{b_1, \ldots, b_\beta\} \cup C)$
and hence $\acl^n(BC) = \acl^n((\{b_1, \ldots, b_\beta\} \cup C)$
it follows that $b_1, \ldots, b_\beta, a_1, \ldots, a_\alpha$ is an $n$-ccs for $A/C$
with core sequence of indices 
$0 = k_0 < k_1 < \ldots < k_s < k_s + l_1 < k_s + l_2 < \ldots k_s + l_t = \beta + \alpha$.
By Lemma~\ref{for every n-cs there is an n-ccs},
we must have $\rk^n(A / C) = \beta + \alpha = \rk^n(B / C) + \rk^n(A / BC)$.
\hfill $\square$

\begin{cor}\label{corollary to modularity}
If $a_1, \ldots, a_k \in M\meq$  and $C \subseteq M\meq$, then 
\begin{align*}
\rk^n(a_1, \ldots, a_k / C) &= \rk^n(a_1 / C) + \rk^n(a_2 / \{a_1\} \cup C) + \ldots + 
\rk^n(a_k / \{a_1, \ldots, a_{k-1}\} \cup C) \\
&\leq \rk^n(a_1 / C) + \ldots + \rk^n(a_k / C).
\end{align*}
\end{cor}

\section{Connection to rosiness}\label{Connection to rosiness}

\noindent
In the previous section we proved that if the algebraic closure satisfies the exchange property on
elements of $M_n$ with $n$-rank 1 (over some $C \subseteq M\meq$), 
then $n$-independence is an independence relation {\em restricted to $M_n$}
(recall Definition~\ref{definition of independence relation}).
In this section we will use this result to make conclusions about (super)rosiness of $\omega$-categorical theories.
In order to do this we also need to involve the notion of thorn-independence.

\begin{defin}\label{definition of thorn-independence}{\rm \cite{EO, Ons}
Let $\bar{a}, \bar{b}$ be finite tuples of elements from $M\meq$ and let $C \subseteq M\meq$.
\begin{enumerate}
\item A formula $\varphi(\bar{x}, \bar{a})$ (with all parameters listed by $\bar{a})$)
{\em strongly divides} over $C$ if $\tp(\bar{a} / C)$ is nonalgebraic (i.e. has infinitely many realizations)
and, for some $k < \omega$, the set of formulas $\{\varphi(\bar{x}, \bar{a}') : \tp(\bar{a}' / C) = \tp(\bar{a} / C) \}$
is $k$-inconsistent (meaning that any set of $k$ of these formulas is inconsistent, in a model of $T\meq$).

\item A formula $\varphi(\bar{x}, \bar{a})$ {\em thorn-divides} over $C$ if there is a finite tuple
$\bar{d}$ of elements from $M\meq$ such that $\varphi(\bar{x}, \bar{a})$
strongly divides over $C\bar{d}$ (= $C \cup \rng(\bar{d})$).

\item A formula $\varphi(\bar{x}, \bar{a})$ {\em thorn-forks} over $C$ if it implies (modulo $T\meq$)
a finite disjunction of formulas, all of which thorn-divide over $C$.

\item For $A \subseteq M\meq$, a complete type $p(\bar{x})$ over $A$ {\em thorn-divides} ({\em thorn-forks}) over $C$
if there is a formula in $p(\bar{x})$ which thorn-divides (thorn-forks) over $C$.

\item We say that $\bar{a}$ is {\em thorn-independent from $\bar{b}$ over $C$}, denoted
$\bar{a} \underset{C}{\thind} \bar{b}$, if $\tp(\bar{a} / C\bar{b})$ does not thorn-fork over $C$.
\end{enumerate}
}\end{defin}

\noindent
The following technical lemma appears in \cite[Remark~3.2]{EO}, but we prove it to make the arguments that
follow self-contained.

\begin{lem}\label{technical lemma about strong dividing}
Let $a, b \in M\meq$, $C \subseteq M\meq$ and let $\bar{c}$ be a sequence of elements from $C$.
If $M\meq \models \varphi(a, b, \bar{c})$ and $\varphi(x, b, \bar{c})$ strongly divides over $C$,
then $b \in \acl\meq(aC) \setminus \acl\meq(C)$.
\end{lem}

\noindent
{\bf Proof.}
The assumption that $\varphi(x, b, \bar{c})$ strongly divides over $C$ means that
$\tp(b / C)$ is nonalgebraic and, for some $k < \omega$,
$\{\varphi(x, b', \bar{c}) : \tp(b' / C) = \tp(b / C)\}$ is $k$-inconsistent.
As explained in \cite[Remark~2.1.2]{Ons}, 
it follows from a compactness argument that there is $\theta(y, \bar{d}) \in \tp(b / C)$
such that 
\[
\Phi = \{\varphi(x, b', \bar{c}) : M\meq \models \theta(b', \bar{d}) \} \ \ \text{is $k$-inconsistent.}
\]
Since $\tp(b / C)$ is nonalgebraic it follows that $b \notin \acl\meq(C)$.
For a contradiction, suppose that $b \notin \acl\meq(aC)$.
Then there are distinct $b_i$, for $i < \omega$, such that 
$\tp(a, b_i, \bar{c}, \bar{d}) = \tp(a, b, \bar{c}, \bar{d})$ for all $i$.
Then $M\meq \models \varphi(a, b_i, \bar{c}) \wedge \theta(b_i, \bar{d})$ for all $i < \omega$,
and this contradicts that $\Phi$ is $k$-inconsistent.
\hfill $\square$

\begin{defin}\label{symmetry restricted to M-n}{\rm
We say that $\indn$ is {\em symmetric restricted to $M_n$}
if for all $a, b \in M_n$,  and $C \subseteq M\meq$, 
if $a \underset{C}{\indn} b$ then $b \underset{C}{\indn} a$.
}\end{defin}

\begin{lem}\label{a technical property of n-independence}
Let $n < \omega$, $a, b \in M_n$, and $C \subseteq M\meq$.
If $b \in \acl^n(aC) \setminus \acl^n(C)$ 
and $\indn$ is symmetric restricted to $M_n$, then $a \underset{C}{\nindn} b$.
\end{lem}

\noindent
{\bf Proof.}
Suppose that $n < \omega$, $a, b \in M_n$, $C \subseteq M\meq$, and
$b \in \acl^n(aC) \setminus \acl^n(C)$.
Then (by the definition of $\rk^n$) $\rk^n(b / C) \geq 1$ and $\rk^n(b / aC) = 0$.
It follows that $b \underset{C}{\nindn} a$, so by the assumed symmetry of $\indn$ restricted to $M_n$
we have $a \underset{C}{\nindn} b$.
\hfill $\square$

\medskip
\noindent
Below we repeat Assumption~\ref{assumption about exchange property}
which implies that $\indn$ is symmetric restricted to $M_n$.

\begin{assump}\label{assumption about exchange property for all n} 
{\bf (Exchange property with respect to $n$)} {\rm
Let $n < \omega$ and suppose that if $C \subseteq M\meq$, $2 \leq k < \omega$, $a_1, \ldots, a_k \in M_n$,
$\rk^n(a_i / C) = 1$ for all $i = 1, \ldots, k$, and
$a_k \in \acl\meq(\{a_1, \ldots, a_{k-1}\} \cup C) \setminus \acl\meq(\{a_2, \ldots, a_{k-1}\} \cup C)$,
then $a_1 \in \acl\meq(\{a_2, \ldots, a_k\} \cup C)$.
}\end{assump}

\noindent
The proof of the next lemma is an adaptation of the proof of Theorem~3.3 in \cite{EO}.

\begin{lem}\label{thorn dependence implies n-dependence}
Suppose that 
Assumption~\ref{assumption about exchange property for all n}
holds for all $n<\omega$. 
Let $\bar{a}, \bar{b} \in M\meq$ and let $C \subseteq M\meq$.
If $\bar{a} \underset{C}{\thnind} \bar{b}$ then $\bar{a} \underset{C}{\nindn} \bar{b}$
for all sufficiently large $n < \omega$.
\end{lem}

\noindent
{\bf Proof.}
Let $m < \omega$, let $\bar{a}$ and
$\bar{b}$ be finite sequences of elements from $M_m$, and let $C \subseteq M\meq$.
First we show that if $\tp(\bar{a} / \bar{b}C)$ thorn divides over $C$, then $\bar{a} \underset{C}{\nindn} \bar{b}$
for all $n \geq m$.
For a contradiction, let $m \leq n < \omega$ and suppose that
$\tp(\bar{a} / \bar{b}C)$ thorn divides over $C$ and $\bar{a} \underset{C}{\indn} \bar{b}$.
Then some $\varphi(\bar{x}, \bar{b}, \bar{c}) \in \tp(\bar{a} / \bar{b}C)$, where $\bar{c} \in C$, 
thorn divides over $C$.
This means that there is $D \supseteq C$ such that $\varphi(\bar{x}, \bar{b}, \bar{c})$
strongly divides over $D$.

The assumptions that $\bar{a} \underset{C}{\indn} \bar{b}$
and that $\rng(\bar{a}) \subseteq M_n$ together with the extension property of $\indn$ 
(Lemma~\ref{extension property of n-independence})
implies that there is $\bar{a}'$ such that $\tp(\bar{a}' / \bar{b}C) = \tp(\bar{a} / \bar{b}C)$ and
$\bar{a}' \underset{\bar{b}C}{\indn} D$, so in particular $\rng(\bar{a}') \subseteq M_n$.
The assumption that  $\bar{a} \underset{C}{\indn} \bar{b}$ (and automorphism invariance of $\indn$)
gives $\bar{a}' \underset{C}{\indn} \bar{b}$.
By transitivity of $\indn$ we then get $\bar{a}' \underset{C}{\indn} \bar{b}D$.
Due to the saturation assumption about $\mcM\meq$ and since $\tp(\bar{a} / \bar{b}C) = \tp(\bar{a}' / \bar{b}C)$ 
there is an elementary function $\sigma$, including $\bar{a}'\bar{b}CD$ in its domain, such that
$\sigma$ fixes $\bar{b}C$ pointwise and $\sigma(\bar{a}') = \bar{a}$.
By automorphism invariance of $\indn$ we get
$\bar{a} \underset{C}{\indn} \bar{b}\sigma(D)$.
As strong dividing is preserved by elementary functions,
$\varphi(\bar{x}, \bar{b}, \bar{c})$ strongly divides over $\sigma(D)$.
Without loss of generality, rename $\sigma(D)$ as $D$.
Then $\bar{a} \underset{C}{\indn} \bar{b}D$ and
$\varphi(\bar{x}, \bar{b}, \bar{c})$ strongly divides over $D$.

From $\bar{a}'\underset{C}{\indn} \bar{b}D$ and monotonicity of $\indn$ we get
$\bar{a} \underset{D}{\indn} \bar{b}D$ which means the same as $\bar{a} \underset{D}{\indn} \bar{b}$.

Since $\mcM\meq \models \varphi(\bar{a}, \bar{b}, \bar{c})$ 
and $\varphi(\bar{x}, \bar{b}, \bar{c})$ strongly divides over $D$ it follows from 
Lemma~\ref{technical lemma about strong dividing}
that $\bar{b} \in \acl\meq(\bar{a}D) \setminus \acl\meq(D)$.
Since $\bar{b} \in M_n$ it follows that 
$\bar{b} \in \acl^n(\bar{a}D) \setminus \acl^n(D)$.
Since Assumption~\ref{assumption about exchange property for all n}
and Proposition~\ref{symmetry of n-independence}
imply that $\indn$ has the symmetry property restricted to $M_n$, it follows from
Lemma~\ref{a technical property of n-independence}
that $\bar{a} \underset{D}{\nindn} \bar{b}$ and this contradicts the conclusion above that 
$\bar{a}  \underset{D}{\indn} \bar{b}$.
Now we have proved that 
\begin{equation}\label{eq dividing implies n-dep}
\text{If $m \leq n$, $\bar{a}, \bar{b} \in M_m$, and $\tp(\bar{a} / \bar{b}C)$ thorn divides over $C$, then $\bar{a} \underset{C}{\nindn} \bar{b}$.}
\end{equation}

Now suppose that $\bar{a}, \bar{b} \in M\meq$, $C \subseteq M\meq$,
and $\bar{a} \underset{C}{\thnind} \bar{b}$, that is,
$\tp(\bar{a} / \bar{b}C)$ thorn forks over $C$.
Then there is $\psi(\bar{x}, \bar{b}, \bar{c})$ in $\tp(\bar{a} / \bar{b}C)$ which thorn forks,
which means that  $\psi(\bar{x}, \bar{b}, \bar{c})$ implies a finite disjunction,
say $\bigvee_{i=1}^s \varphi_i(\bar{x}, \bar{b}, \bar{c}, \bar{d}_i)$ of formulas, all of which 
thorn divide over $C$. 
Let $D$ be the union of all $\rng(\bar{d}_i)$.
Then every extension of $\tp(\bar{a} /\bar{b}C)$ to a complete type $p(\bar{x})$ over $\bar{b}CD$
thorn divides over $C$. 
As $D$ is finite there is $m < \omega$ such that all elements in $\bar{a}\bar{b}D$ belong to $M_m$.
Let $n \geq m$.
By~(\ref{eq dividing implies n-dep}),
whenever $\tp(\bar{a}' / \bar{b}CD)$ extends $\tp(\bar{a} /\bar{b}C)$ 
then $\bar{a}' \underset{C}{\nindn} \bar{b}D$.
By the extension property of $\indn$ there is $\bar{a}'$ such that 
$\tp(\bar{a}' / \bar{b}C) = \tp(\bar{a} / \bar{b}C)$ and 
$\bar{a}' \underset{\bar{b}C}{\indn} D$.
Now transitivity of $\indn$ implies that $\bar{a}' \underset{C}{\nindn} \bar{b}$.
As $\tp(\bar{a}' / \bar{b}C) = \tp(\bar{a} / \bar{b}C)$ we get  $\bar{a} \underset{C}{\nindn} \bar{b}$.
\hfill $\square$

\begin{prop}\label{thorn independence has local character}
If $T$ is $\omega$-categorical and Assumption~\ref{assumption about exchange property for all n}
holds for all $n < \omega$, then thorn-independence has local character.
\end{prop}

\noindent
{\bf Proof.}
By Remark~\ref{remark on locality} it suffices to prove that 
there do {\em not} exist $a \in M\meq$ and $b_\alpha \in M\meq$, for $\alpha < \aleph_1$, such that 
$a \underset{(b_i)_{i < \alpha}}{\thnind} b_\alpha$ for all $\alpha < \aleph_1$.
Towards a contradiction, suppose that there are $a \in M\meq$ and $b_\alpha \in M\meq$, for $\alpha < \aleph_1$, 
such that $a \underset{(b_i)_{i < \alpha}}{\thnind} b_\alpha$ for all $\alpha < \aleph_1$.
Let $B = \{b_\alpha : \alpha < \aleph_1\}$.
By Lemma~\ref{locality property again},
there is a countable $C \subseteq B$ such that 
$a \underset{C}{\indn} B$ for all $n < \omega$.
Let $\beta = \sup\{\alpha : b_\alpha \in C\}$.
As $\aleph_1$ is a regular cardinal and $C$ is countable it follows that $\beta < \aleph_1$, so $\beta$ is a
countable ordinal.
Let $B_\beta = \{b_\alpha : \alpha \leq \beta\}$ so $B_\beta$ is countable and $C \subseteq B_\beta$.
By the monotonicity of $\indn$ we get $a \underset{B_\beta}{\indn} B$ for all $n < \omega$.
By assumption we have $a \underset{B_\beta}{\thnind} b_{\beta + 1}$, 
so Lemma~\ref{thorn dependence implies n-dependence}
(which uses Assumption~\ref{assumption about exchange property for all n})
implies that $a \underset{B_\beta}{\nindn} b_{\beta + 1}$ for all sufficiently large $n$.
Since $b_{\beta + 1} \in B$ it follows that $a \underset{B_\beta}{\nindn} B$ for all sufficiently large $n$.
But this contradicts the earlier conclusion that  $a \underset{B_\beta}{\indn} B$ for all $n < \omega$.
\hfill $\square$

\medskip
\noindent
The following result by Ealy and Onshuus  \cite{EO} is crucial for the results that follow:

\begin{theor}\label{local character of thorn forking implies rosiness} {\rm \cite[Theorem 3.7]{EO}}
A theory is rosy if and only if thorn indepence has local character.
\end{theor}

\begin{theor}\label{the assumption implies that T is rosy}
If $T$ is $\omega$-categorical and Assumption~\ref{assumption about exchange property for all n}
holds for all $n < \omega$,
then $T$ is rosy.
\end{theor}

\noindent
{\bf Proof.}
If  $T$ is $\omega$-categorical and Assumption~\ref{assumption about exchange property for all n}
holds for all $n < \omega$, then, by Proposition~\ref{thorn independence has local character},
thorn independence has local character.
By Theorem~\ref{local character of thorn forking implies rosiness},
$T$ is rosy.
\hfill $\square$

\begin{defin}\label{definition of trivial acl up to n}{\rm
(i) We say that algebraic closure is {\em trivial in $T\meq$} if whenever $\mcM \models T$,
$a \in M\meq$, $B \subseteq M\meq$ and $a \in \acl\meq(B)$, then $a \in \acl\meq(b)$ for some $b \in B$.\\
(ii) Let $n < \omega$. We say that algebraic closure in $T\meq$ is {\em trivial up to  level $n$}
if the following holds for every $\mcM \models T$:
If $a \in M_n$, $B \subseteq M_n$, $C \subseteq M\meq$, and $a \in \acl\meq(BC) \setminus \acl\meq(C)$
then $a \in \acl\meq(bC)$ for some $b \in B$.
}\end{defin}

\begin{lem}\label{trivial acl up to m}
Let $T$ be $\omega$-categorical and $m < \omega$.
If algebraic closure in $T\meq$ is trivial up to level $m$ then 
Assumption~\ref{assumption about exchange property for all n} holds for all $n \leq m$.
\end{lem}

\noindent
{\bf Proof.}
Suppose that $T$ is $\omega$-categorical and that algebraic closure in $T\meq$ is trivial up to some $m < \omega$.
Fix any $n \leq m$ and suppose that
$C \subseteq M\meq$, $2 \leq k < \omega$, $a_1, \ldots, a_k \in M_n$,
$\rk^n(a_i / C) = 1$ for all $i = 1, \ldots, k$, and
$a_k \in \acl\meq(\{a_1, \ldots, a_{k-1}\} \cup C) \setminus \acl\meq(\{a_2, \ldots, a_{k-1}\} \cup C)$.
Since $M_n \subseteq M_m$ and 
the algebraic closure in $T\meq$ is trivial up to $m$ it follows that $a_k \in \acl\meq(a_1)$.

For a contradiction, suppose that $a_1 \notin \acl\meq(\{a_k\} \cup C)$.
Then $\rk^n(a_1 / \{a_k\} \cup C) \geq 1$.
Since $\rk^n(a_k / C) = 1$ (and $a_k \in M_n$) we have $a_k \notin \acl\meq(C)$.
It now follows from the definition of $\rk^n$ that $\rk^n(a_1 / C) \geq 2$ which contradicts the assumption.
Hence $a_1 \in \acl\meq(\{a_k\} \cup C)$, which since $a_1 \notin \acl\meq(C)$
(and algebraic closure in $T\meq$ is trivial up to $m$)
implies that $a_1 \in \acl\meq(a_k)$.
\hfill $\square$

\begin{theor}\label{trivial acl implies that T is rosy}
If $T$ is $\omega$-categorical with trivial algebraic closure in $T\meq$
then $T$ is rosy.
\end{theor}

\noindent
{\bf Proof.}
Suppose that $T$ is $\omega$-categorical with trivial algebraic closure in $T\meq$.
Then, for every $n < \omega$, algebraic closure in $T\meq$ is trivial up to $n$.
Then Lemma~\ref{trivial acl up to m} implies that 
Assumption~\ref{assumption about exchange property for all n} 
holds for all $n < \omega$.
Now Theorem~\ref{the assumption implies that T is rosy}
implies that $T$ is rosy.
\hfill $\square$

\begin{defin}\label{definition of thorn rank}{\rm \cite{EO, Ons} 
Let $p(\bar{x})$ be a complete type over some $A \subseteq M\meq$
(and we assume that consistency is part of the definition of being a type).
\begin{enumerate}
\item We define the $U^{\text{\rm \th}}$-rank of $p(\bar{x})$ as follows:
\begin{enumerate}
\item $U^{\text{\rm \th}}(p(\bar{x}) \geq 0$.

\item For every ordinal $\alpha$, $U^{\text{\rm \th}}(p(\bar{x})) \geq \alpha + 1$
if there is $b \in M\meq$ and a complete type $q(\bar{x})$ over $Ab$ such that
$p(\bar{x}) \subseteq q(\bar{x})$,
$U^{\text{\rm \th}}(q(\bar{x})) \geq \alpha$, 
and $q(\bar{x})$ thorn forks over $A$.

\item For every limit ordinal $\beta$, $U^{\text{\rm \th}}(p(\bar{x})) \geq \beta$
if  $U^{\text{\rm \th}}(p(\bar{x})) \geq \alpha$ for all ordinals $\alpha < \beta$.

\item For every ordinal $\alpha$,  $U^{\text{\rm \th}}(p(\bar{x})) = \alpha$ if 
 $U^{\text{\rm \th}}(p(\bar{x})) \geq \alpha$ and  $U^{\text{\rm \th}}(p(\bar{x})) \not\geq \alpha + 1$.
\end{enumerate}
\item A theory is {\em superrosy} if for every complete type $p(x)$ over $\es$
there is an ordinal $\alpha$ such that  $U^{\text{\rm \th}}(p(x)) = \alpha$.
If in addition  $U^{\text{\rm \th}}(p(x)) < \omega$ for every complete type $p(x)$
over $\es$, then we say that the theory has {\em finite $U^{\text{\rm \th}}$-rank}.
\end{enumerate}
}\end{defin}

\noindent
Since thorn forking has the extension property and since
$U^{\text{\rm \th}}(\bar{a} / B) \leq U^{\text{\rm \th}}(\bar{a} / A)$
if $A \subseteq B$ (see e.g. \cite[Theorem 2.7]{EO})
it follows from \cite[Fact ~4.4]{EO} that the definition of superrosy theory above is equivalent
to the one given in \cite[Section~4.1]{EO}.

\begin{theor}\label{concluding that T is superrosy}
If $T$ is $\omega$-categorical with geometric elimination of imaginaries and 
Assumption~\ref{assumption about exchange property for all n}
holds for $n = 0$,
then $T$ is superrosy with finite $U^{\text{\rm \th}}$-rank.
\end{theor}

\noindent
{\bf Proof.}
Let $T$ be $\omega$-categorical with geometric elimination of imaginaries and suppose that 
Assumption~\ref{assumption about exchange property for all n}
holds for $n = 0$.
It follows from Theorem~\ref{n-independence is an independence relation}
that $\indo$ is an independence relation.
Since $T$ has geometric elimination of imaginaries
Lemma~\ref{n-independence and geometric elimination of imaginaries}
implies that, for all $n < \omega$ and all $A, B, C \subseteq M\meq$, 
$A \underset{C}{\indn} B$ if and only if $A \underset{C}{\indo} B$.
Hence, for all $n < \omega$, $\indn$ is an independence relation.
It follows that 
Assumption~\ref{assumption about exchange property for all n} 
holds for all $n < \omega$, so 
Lemma~\ref{thorn dependence implies n-dependence}
implies that if $a \underset{C}{\thnind} b$ then $a \underset{C}{\nindn} b$ for all sufficiently large $n$.

To show that $T$ is superrosy with finite $U^{\text{\rm \th}}$-rank it suffices to show that 
$U^{\text{\rm \th}}(p) < \omega$ for every complete type $p(x)$ over $\es$.
For a contradiction, suppose that $U^{\text{\rm \th}}(p) \geq \omega$.
Let $\alpha = \rk^0(a)$ where $a$ is any realization of $p(x)$, so $\alpha < \omega$.
Take any $\alpha < \beta < \omega$. As $U^{\text{\rm \th}}(p) \geq \beta$ there are $a \in M\meq$ realizing $p(x)$ and
$b_k \in M\meq$, for $k \leq \beta$,
such that
$a \underset{(b_i)_{i < k}}{\thnind} b_k$ for all $k \leq \beta$.
From our conclusion above it follows that there is $n' < \omega$ such that 
if $n \geq n'$ then $a \underset{(b_i)_{i < k}}{\nindn} b_k$ for all $k \leq \beta$.
As $\indn$ coincides with $\indo$ it follows that 
$a \underset{(b_i)_{i < k}}{\nindo} b_k$ for all $k \leq \beta$.
Hence 
\[
\rk^0(a) > \rk^0(a / b_0) > \rk^0(a / b_0, b_1) > \ldots > \rk^0(a / b_0, \ldots, b_\beta) \geq 0
\]
so $\rk^0(a) \geq \beta > \alpha$, contradicting the choice of $\alpha$.
\hfill $\square$

\begin{rem}\label{remark on concluding that T is superrosy}{\rm
Due to Lemma~\ref{trivial acl up to m}, the conclusion of 
Theorem~\ref{concluding that T is superrosy} still holds if the assumption
that Assumption~\ref{assumption about exchange property for all n}
holds for $n = 0$
is replaced by the assumption that algebraic closure in $T\meq$ is trivial up to level 0.
}\end{rem}

\subsection*{Epilogue}

The work resulting in this article began by considering a sequence $(\mcB_n : n < \omega)$ of finite structures $\mcB_n$
where $\lim_{n\to\infty}|B_n| = \infty$ and where all $\mcB_n$ have a ``uniformly well-behaved'' closure operator,
where ``well-behaved'' essentially means that the closure operator is uniformly definable and
has the properties of the algebraic closure
operator in an $\omega$-categorical structure. The idea was that such a sequence $(\mcB_n : n \in \mbbN)$
would generalize the ``base sequence'' of finite trees with bounded height considered in \cite{KT} where a closure operator
 has a crucial role (the operator that collects all ancestors of a set of vertices of the tree).
It turns out that if we have such a well-behaved closure operator and the sequence $(\mcB_n : n \in \mbbN)$
has an additional property (not explained here), then the same sequence has an infinite (countable) limit structure
(in a seemingly new sense), say $\mcM$,
which is $\omega$-categorical. 
However, only knowing that the limit $\mcM$ is $\omega$-categorical did not seem 
to be of much interest. It was clear that the closure operator on the structures $\mcB_n$ gives rise to a 
definable closure operator $\cl$ on $\mcM$ which is such that $\cl(A) \subseteq \acl(A)$ for all $A \subseteq M$
where $\acl(A)$ denotes the algebraic closure of $A$ in $\mcM$.
Since the ``nice'' properties of $\cl$ on all $\mcB_n$ meant that a notion of rank was possible to define
on all $\mcB_n$, in a uniform way, it followed that a notion of rank (defined in terms of $\cl$) existed in $\mcM$.
However, if the notion of rank is going to be useful to define a ``nice'' independence relation
in the sense of \cite{KP},
then one has to consider ``closure'' not only on ``real elements'' of $\mcM$ but also on imaginary elements.
Since, for all $A \subseteq M$, $\cl(A) \subseteq \acl(A)$ it seemed to be reasonable,
at the present state of affairs at that time,
to just forget about the sequence
$(\mcB_n : n \in \mbbN)$ 
and consider ranks and independence relations defined by the algebraic
closure operator in an $\omega$-categorical structure.
This eventually led to the results presented here.
Given these results it may now be more meaningful to go back to considerations of the
sequence $(\mcB_n : n < \omega)$ and the connections between the closure operator on all $\mcB_n$
(or some generalization of it to take into account ``imaginary elements in $\mcB_n$'')
and the limit structure $\mcM$.

\subsection*{Acknowledgement}

This work is partially
supported by the Swedish Research Council, grant 2023-05238.


\begin{thebibliography}{99}\label{References}

\bibitem{Adl} H. Adler, A geometric introduction to forking and thorn-forking,
{\em Journal of Mathematical Logic}, Vol. 9 (2009) 1--20.

\bibitem{AZ} G. Ahlbrandt, M. Ziegler, Quasi-finitely axiomatizable totally categorical theories,
{\em Annals of Pure and Applied Logic}, Vol. 30 (1986), 63-82.

\bibitem{BFM} J. Baldwin, J. Freitag, S. Mutchnik,
Simple homogeneous structures and indiscernible sequence invariants,
\url{https://arxiv.org/abs/2405.08211}.

\bibitem{CHL} G. Cherlin, L. Harrington, A. H. Lachlan, $\omega$-categorical $\omega$-stable structures,
{\em Annals of Pure and Applied Mathematics}, Vol. 28 (1986) 103--135.

\bibitem{CH} G. Cherlin, E. Hrushovski, Finite Structures with Few Types, 
{\em Annals of Mathematics Studies}, Nr. 152, Princeton University Press (2003).

\bibitem{Con} G. Conant, An axiomatic approach to free amalgamation,
{\em The Journal of Symbolic Logic}
Vol. 82, (2017), 648--671.

\bibitem{EO} C. Ealy, A. Onshuus, Characterizing rosy theories, 
{\em The Journal of Symbolic Logic}, Vol. 72 (2007) 919--940.

\bibitem{Eva} D. Evans, $\aleph_0$-categorical structures with a predimension,
{\em Annals of Pure and Applied Logic}, Vol. 116 (2002) 157--186.

\bibitem{Gag} J. Gagelman, Stability in geometric theories,
{Annals of Pure and Applied Logic}, Vol. 132 (2005) 313--326.

\bibitem{Hod} W. Hodges,  {\em Model theory}, Cambridge University Press (1993).

\bibitem{HHM} W. Hodges, I.M. Hodkinson, D. Macpherson, Omega-categoricity, relative categoricity
and coordinatisation, 
{\em Annals of Pure and Applied Logic}, Vol. 46 (1990) 169--199.

\bibitem{Hru1} E. Hrushovski, Totally categorical structures
{\em Transactions of the American Mathematical Society},
Vol. 313, (1989), 131--159.

\bibitem{Hru2} E. Hrushovski, An new strongly minimal set,
{\em Annals of Pure and Applied Logic}, Vol. 62 (1993) 147--166.

\bibitem{KLM} W. M. Kantor, M. W. Liebeck, H. D. Macpherson, 
$\aleph_0$-categorical structures smoothly approximated by finite substructures,
{\em Proceedings of the London Mathematical Society}, Vol. 59 (1989) 439--463.

\bibitem{Kim_book} B. Kim, {\em Simplicity Theory}, Oxford University Press (2014).

\bibitem{KP} B. Kim, A. Pillay, Simple theories, 
{\em Annals of Pure and Applied Logic}, Vol. 88 (1997) 149--164.

\bibitem{Kop} V. Koponen, Binary simple homogeneous structures,
{\em Annals of Pure and Applied Logic},
Vol. 169 (2018) 1335--1368.

\bibitem{KT} V. Koponen, Y. Tousinejad, Random expansions of trees with bounded height,
{\em Theoretical Computer Science},
Vol. 1040 (2025) 115201.

\bibitem{Lach97} A. H. Lachlan, Stable finitely homogeneous structures: a survey, 
in B. T. Hart, A. H. Lachlan, M. A. Valeriote (editors), 
{\em Algebraic Model Theory}, Kluwer Academic Publishers (1997) 145--159.

\bibitem{Mor} M. Morley, Categoricity in power,
{\em Transactions of the American Mathematical Society}, Vol. 114 (1965), 514--538.

\bibitem{NW} L. Newelski, R. Wencel, Definable sets in Boolean-ordered o-minimal structures. I, 
{\em The Journal of Symbolic Logic}, Vol. 66 (2001) 1821--1836.


\bibitem{Ons} A. Onshuus, Properties and consequences of thorn-independence,
{\em The Journal of Symbolic Logic}, Vol. 71 (2006) 1--21.

\bibitem{OS} A. Onshuus, P. Simon, Dependent finitely homogeneous rosy theories,
\url{https://arxiv.org/pdf/2107.02727}.

\bibitem{Ox} J. Oxley, {\em Matroid Theory (2nd edn)}, 
Oxford University Press (2011).

\bibitem{She} S. Shelah, {\em Classification Theory}, Revised Edition, North-Holland (1990).

\bibitem{Sim_book} P. Simon, {\em A guide to NIP theories}, 
Lecture Notes in Logic Vol. 44, Cambridge University Press (2015).

\bibitem{Sim22} P. Simon, NIP $\omega$-categorical structures: The rank 1 case,
{\em Proceedings of the London Mathematical Society}, 
Vol. 125 (2022) 1253--1331.

\bibitem{TZ} K. Tent, M. Ziegler, {\em A Course in Model Theory},
Lecture Notes in Logic, Vol. 40, Cambridge University Press (2012).

\bibitem{Van} L. van den Dries, Tame Topology and O-minimal Structures,
{\em London Mathematical Society Lecture Note Series}, Vol. 248, 
Cambridge University Press (1998).

\bibitem{Zil} B. Zilber, {\em Uncountably Categorical Theories}, AMS Translations of Mathematical Monographs,
Vol. 117 (1993).


\end{thebibliography}
\end{document}